\documentclass[11pt,a4paper]{article}
\usepackage[utf8]{inputenc}
\usepackage{amsmath,bm}
\usepackage{enumerate}
\usepackage{amsmath,nccmath}
\usepackage{amsthm}
\usepackage{amsfonts}
\usepackage[square, comma, sort&compress,numbers]{natbib}
\makeatletter  
\renewcommand{\@biblabel}[1]{[#1]\hfill}
\makeatother
\usepackage{amssymb}
\usepackage{setspace}
\usepackage{graphicx}
\usepackage{epstopdf}
\usepackage[shortlabels]{enumitem} 
\usepackage{verbatim}
\usepackage{color}
\usepackage[backref=page]{hyperref}
\hypersetup{
     colorlinks=true,    
     linktocpage=true,
    }
\renewcommand*{\backref}[1]{}
\renewcommand*{\backrefalt}[4]{~~{\tiny%
    \ifcase #1 Not cited.%
          \or [Cited on page~#2.]%
          \else [Cited on pages #2.]%
    \fi%
    }}
\usepackage[marginal]{footmisc}
\usepackage[section]{placeins}
\usepackage{indentfirst}
\usepackage[left=2cm,right=2cm,top=2cm,bottom=2cm]{geometry}
\usepackage{float}
\usepackage[doipre={DOI:~}]{uri} 
\usepackage{enumitem}
\usepackage{tocloft}
\usepackage{titlesec} 
\titleformat{\subsubsection}[runin]
  {\normalfont\normalsize\bfseries}{\thesubsubsection}{1em}{}

\author{}
\theoremstyle{plain}
\newtheorem{prop}{Proposition}[section]
\newtheorem*{prop*}{Proposition}
\newtheorem{thm}{Theorem}
\newtheorem{cor}[thm]{Corollary}
\newtheorem{lem}[prop]{Lemma}
\newtheorem*{lem*}{Lemma}
\newtheorem*{thm*}{Theorem}

\newtheorem*{claim*}{Claim}

\newtheorem{thms}[prop]{Theorem}

\newtheorem*{ques*}{Question}

\theoremstyle{definition}
\newtheorem{rem}[prop]{Remark}
\newtheorem{defi}[prop]{Definition}
\newtheorem*{defi*}{Definition}
\newtheorem{ex}[prop]{Example}

\newtheorem*{conj*}{Conjecture}

\newtheorem*{conv}{Convention}

\allowdisplaybreaks
\usepackage{mathtools}

\DeclarePairedDelimiter\floor{\lfloor}{\rfloor}

\usepackage{authblk}

\renewcommand{\geq}{\geqslant}
\renewcommand{\leq}{\leqslant}

\definecolor{darkred}{rgb}{.56,0,0}

\makeatletter
\newcommand*{\rom}[1]{\expandafter\@slowromancap\romannumeral #1@}
\newcommand\primeitem{%
 \item[(\roman{enumi}\textquotesingle)]\def\@currentlabel{(\roman{enumi}\textquotesingle)}}
\newcommand\pprimeitem{%
 \item[\arabic{enumi}\textquotesingle]\def\@currentlabel{\arabic{enumi}\textquotesingle}}

\newcommand{\setword}[2]{%
  \phantomsection
  #1\def\@currentlabel{\unexpanded{#1}}\label{#2}%
} 
 
\makeatother

\title{{\bf  $C^1$-robust homoclinic tangencies}\\
}
\author{Dongchen Li\footnote{dongchenli@fudan.edu.cn, Shanghai Center for Mathematical Sciences, Fudan University, China}} 

\begin{document}
\def\D{\mathrm{D}}
\def\d{\mathrm{d}}
\def\tr{\mathrm{tr}\,}
\def\diff{\mathrm{Diff}}
\def\eps{\varepsilon}
\def\pp{\varphi}
\def\card{\card}
\def\homo{\overset{\mathrm{h.r.}}{\sim}}
\def\usim{\overset{\mathrm{u}}{\sim}}
\def\susim{\overset{\mathrm{su}}{\sim}}
\def\xsim{\overset{\times}{\sim}}
\def\inte{\operatorname{int}}
\def\cl{\operatorname{cl}}
\def\eff{\mathrm{eff}}
\def\dom{\mathrm{dom}}
\def\B{\mathrm{B}}
\def\C{\mathrm{C}}
\def\L{\mathrm{L}}
\def\R{\mathrm{R}}
\def\G{\mathrm{G}}
\def\s{s}
\def\spc{\quad\mbox{and}\quad}
\def\X{\mathbb{X}}
\def\Y{\mathbb{Y}}
\def\Z{\mathbb{Z}}

\def\ind{\operatorname{ind}}
\def\arg{\operatorname{Arg}}
\def\rank{\operatorname{rank}}
\def\card{\operatorname{card}}
\def\diam{\operatorname{diam}}

\numberwithin{figure}{section}

\maketitle

\par{}
\noindent{\bf Abstract.} 
The aim of this paper is twofold. First, 
we introduce standard blenders (special hyperbolic sets) and their variations,  and prove their fundamental properties on the generation of $C^1$-robust tangencies. In particular, these blenders  appear after $C^r$-small perturbations of any diffeomorphism having a heterodimensional cycle of coindex 1.
Next, as an application, we show that
unfolding a homoclinic tangency  to a hyperbolic periodic point can produce uncountably many $C^1$-robust homoclinic tangencies, provided that either this point  is involved in a coindex-1 heterodimensional cycle, or the central dynamics near it is not essentially two-dimensional. The result answers a question posed by Bonatti and D{\'i}az in \citep{BonDia:12b}. 

\par{}
\noindent {\bf Keywords.} blender, heterodimensional cycle, homoclinic tangency, nonhyperbolic dynamics.
\par{}
\noindent {\bf AMS subject classification.} 37C05, 37C29, 37C75.

\tableofcontents

\section{Introduction}


A {\em homoclinic tangency}  refers to a non-transverse intersection between the stable and unstable manifolds of a saddle\footnote{By {saddles} we mean hyperbolic periodic points with multipliers  inside and also outside the unit  circle (in particular, they are not sinks or sources). The multipliers are the eigenvalues of the linearization matrix of $f^{per(O)}$ at $O$.
}
 $O$ of a diffeomorphism $f$.
  Since the Kupka-Smale Theorem asserts that, $C^r$-generically, the invariant manifolds of saddles intersect transversely, the homoclinic tangency to $O$ cannot be robust.
   However, this situation changes when $O$ is replaced by a non-trivial  hyperbolic basic (i.e, compact, transitive and locally maximal)  set.
\begin{defi}[Robust homoclinic tangencies]\label{defi:homo}
Let $f$ be a $C^r$ diffeomorphism with a non-trivial hyperbolic basic   set $\Lambda$. We say that  $\Lambda$ has a  {\em $C^r$-robust homoclinic tangency}  if there exists a $C^r$-neighborhood $\mathcal{U}$ of $f$ such that for  every $g\in\mathcal{U}$ the continuation $\Lambda_g$ of $\Lambda$ has an orbit of non-transverse intersection between $W^u(\Lambda_g)$ and $ W^s(\Lambda_g)$. 
\end{defi}

%

The first  example of robust homoclinic tangencies was given by Newhouse \citep{New:70}, where he showed the existence of $C^2$-robust homoclinic tangencies for   surface diffeomorphisms. Later, he  proved \citep{New:79} that such robust phenomenon is, in fact,  abundant: 

\noindent {\em every $C^r$ $(r=2,\dots,\infty,\omega)$ surface diffeomorphism which has a homoclinic tangency to a saddle can be approximated in the $C^r$ topology by a diffeomorphism which has a hyperbolic basic  sets exhibiting $C^2$-robust homoclinic tangencies.}

\noindent A generalization to higher dimensions with the same regularity was done  in  \citep{GonTurShi:93c,PalVia:94,Rom:95}. Those open sets of homoclinic tangencies   carry many important dynamical phenomena: the coexistence of infinitely many sinks\footnote{Periodic points with all multipliers lying inside the unit circle.} 
 \citep{New:74}, the super-exponential growth of periodic orbits \citep{Kal:00},  the abundance of strange attractors \citep{MorVia:93}, and the existence of wandering Fatou components in $\mathbb{C}^2$ \citep{BerBie:22}, among others.

It is natural to ask whether the robust homoclinic tangencies in the Newhouse construction are $C^1$-robust.  A key ingredient in Newhouse's proof is  the stable intersection between two Cantor sets obtained from the stable and unstable laminations of two horseshoes with large {\em  thickness}. Roughly speaking, large thickness implies a non-empty intersection between the Cantor sets, and each intersection point corresponds to a heteroclinic tangency between the two horseshoes. The $C^2$-robust tangencies follow from the continuity dependence of thickness on $C^2$  diffeomorphisms. However, Ures showed \citep{Ure:95} that the thickness loses continuous dependence on diffeomorphisms of class  $C^1$, and Moreira further built in \citep{Mor:11} a general theory that no horseshoes (hyperbolic basic sets) of $C^1$ surface diffeomorphisms can have $C^1$-robust tangencies. Since all higher dimensional generalization of the Newhouse theorem more or less rely on the original two dimensional construction, the above-mentioned results strongly indicate that Newhouse's thick horseshoe construction cannot produce $C^1$-robust  tangencies. In this paper, we prove

\begin{thm}\label{thm:main}
Every $C^r$ $(r=1,\dots,\infty,\omega)$ diffeomorphism $f$ which  has a homoclinic tangency to a saddle of effective dimension larger than one can be approximated in the $C^r$ topology by a diffeomorphism which has  a hyperbolic basic set  exhibiting uncountably many $C^1$-robust homoclinic tangencies.
\end{thm}

\begin{rem}
The result shows  the abundance of $C^1$-robust homoclinic tangencies in the {\em $C^r$ topology}, in contrast to the earlier work \citep{BonDia:12b} by Bonatti and D\'iaz on the {\em $C^1$ abundance}, where they proved that, if a diffeomorphism has a  homoclinic class satisfying the index-variation and non-domination conditions, then it can be {\em $C^1$-approximated} by a diffeomorphism with a $C^1$ robust homoclinic tangency. 
\end{rem}

This theorem follows directly from a more detailed version -- Theorem~\ref{thm:ht} in Section~\ref{sec:introNewhouse}. 
The notion of effective dimension, denoted by $d_{\eff}$, will be made precise in Definition \ref{defi:d_eff}, and the condition  $d_{\eff}>1$   basically means that  the central dynamics near the saddle is not essentially two-dimensional.  In particular, Theorem~\ref{thm:main}  answers a question of Bonatti and D\'iaz  concerning the sufficient conditions for creating $C^1$-robust homoclinic tangencies inside the Newhouse domain (see \citep[Question 1.11]{BonDia:12b}).

As mentioned, an important consequence of the Newhouse theorem is the prevalence of the phenomenon of  infinitely many coexisting sinks: $C^r$-close to any diffeomorphism having a homoclinic tangency to a saddle that can produce sinks, there is a $C^2$-open set $\mathcal{U}$  where $C^r$-generic maps display  infinitely many sinks.  Theorem~\ref{thm:main} optimizes this result in the sense that, when $d_{\eff}>1$, the set $\mathcal{U}$ can be taken $C^1$-open:
\begin{cor}\label{cor:main}
If the saddle in Theorem~\ref{thm:main} is of sink-producing type, then the diffeomorphism $f$ lies in the $C^r$ closure  of a $C^1$-open set $\mathcal{U}$ such that $\mathcal{U}$ has a  $C^r$-residual subset $\mathcal{R}$ where every diffeomorphism displays infinitely many sinks.
\end{cor}

The corollary is proved in Section~\ref{sec:introNewhouse}.
Saddles of sink-producing type were introduced by Gonchenko, Shilnikov and Turaev  in \citep{GonShiTur:08} and will be specified in Definition~\ref{defi:sink}. Roughly speaking, a saddle is of sink-producing type if it has no strong-unstable multipliers and $f$ contracts volumes in the central subspace, which in particular includes the sectionally-dissipative case considered in \citep{New:74, GonTurShi:93c, PalVia:94}.

\subsection{Robustness of tangencies: from $C^2$  to $C^1$}

In contrast to the essentially two-dimensional Newhouse construction, there are different constructions in dimension three or higher that lead to $C^1$-robust homoclinic tangencies. Examples with specific global constructions were presented in \citep{Sim:72b,Asa:08}. A more general local approach was proposed by Bonatti and D\'iaz   in \citep{BonDia:12b} using the  so-called  {\em blender-horseshoes}, which are a special type of  hyperbolic basic sets and play a role analogous to the thick horseshoes in the Newhouse construction.  They proved that a homoclinic tangency to the blender-horseshoe can be made $C^1$-robust by an arbitrarily $C^r$-small perturbation. 

Inspired by the work \citep{BonDia:12b}, we introduce the notions of {\em separated standard blenders} and {\em arrayed standard blenders}, which are variations of the {\em standard blenders} defined  in \citep{LiTur:24}; precise definitions are given in Sections \ref{sec:sepa} and \ref{sec:array}, respectively. For brevity, we will also refer to them as separated and arrayed blenders. They are both given as the locally maximal invariant sets induced from  partially hyperbolic Markov partitions (see Definition \ref{defi:markov}). In particular, separated  blenders are a  generalization of blender-horseshoes. 
  Let $\mathcal{M}$ be a manifold of dimension at least three and denote by $\diff^r(\mathcal{M})$ the space of $C^r$ diffeomorphisms of $\mathcal{M}$.  
\begin{thm}\label{thm:tangency}
Let $f\in \diff^r(\mathcal{M})$, $r=1,\dots,\infty,\omega$,  have a separated standard blender $\Lambda$ with a  homoclinic tangency. Then, $f$ can be approximated in the $C^r$ topology by a diffeomorphism $g$ such that $\Lambda_g$ exhibits a $C^1$-robust homoclinic tangency.
\end{thm}

We prove this theorem  in Section~\ref{sec:equi}.  The arrayed  blenders carry a more delicate structure and can produce a large number of robust tangencies. 
\begin{thm}\label{thm:uncount}
Let $f\in \diff^r(\mathcal{M})$, $r=1,\dots,\infty,\omega$,  have  an arrayed  standard blender $\Lambda$ with a  homoclinic tangency. Then, $f$ can be approximated in the $C^r$ topology by a diffeomorphism $g$ such that $\Lambda_g$ exhibits uncountably many $C^1$-robust homoclinic tangencies.
\end{thm}



This result is  proved in   Section~\ref{sec:robuncon}. It is important to know that separated/arrayed blenders   can be produced through $C^r$-small perturbations (see Section~\ref{sec:hthdc}), whereas the creation of blender-horseshoes, in all known works, requires $C^1$ perturbations that are large with respect to the $C^r$ topology. We provide in Section~\ref{sec:differ}  a brief explanation of the mechanism for generating robust tangencies and a comparison among blender-horseshoes, separated blenders and arrayed blenders.


\subsection{Creating $C^1$-robust homoclinic tangencies from heterodimensional cycles}\label{sec:hthdc}

We now state the results on the creation of standard blenders from   heterodimensional cycles. 
Let us recall some definitions. The index of a transitive hyperbolic set, denoted by $\ind(\cdot)$, refers to the rank of its unstable bundle. 

\begin{defi}[Heterodimensional cycles]\label{defi:hdc}
We say that a diffeomorphism has a {\em heterodimensional cycle} involving two saddles $O_1$ and $O_2$ if $\ind(O_1)\neq \ind(O_2)$, and  $W^u(O_1)\cap W^s(O_2)\neq\emptyset $ and $W^u(O_2)\cap W^s(O_1)\neq\emptyset $. The difference of the indices is called the {\em coindex} of the cycle.
\end{defi}

 We enumerate the involved saddles in the order   that  $\ind(O_1)+1= \ind(O_2)$.
The {\em center-stable} and  {\em center-unstable} multipliers of $O_i$ $(i=1,2)$ are those nearest to the unit circle from inside and, respectively, from outside. In general position, 
each $O_i$ has either  one simple and real center-stable multiplier or two, forming a conjugate complex (nonreal) pair. Similar for the center-unstable multiplier.   
Denote by $\lambda$ the center-stable multiplier of $O_1$ and by $\gamma$ the center-unstable multiplier of $O_2$. Then,  a heterodimensional cycle  falls into one of the three different cases:
\begin{itemize}[nosep]
\item {\em saddle}: both $\lambda$ and $\gamma$ are real;
\item {\em saddle-focus}: one of $\lambda$ and $\gamma$ is real and the other is not; and
\item {\em double-focus}: both $\lambda$ and $\gamma$ are not real.
\end{itemize}

Standard blenders  carry a partially hyperbolic structure with a one-dimensional center. They are further called {\em center-stable (cs)} or {\em center-unstable (cu)}, depending on whether the central dynamics is contracting or expanding; see Definition \ref{defi:blender}. The standard blenders near a heterodimensional cycle are obtained as invariant sets of the first return map along the cycle. In the saddle-focus and double-focus cases, there is flexibility  (due to the rotations brought by complex central multipliers) for realizing both contracting and expanding central dynamics. However, in the saddle case, this  contraction/expansion is relatively rigid and is governed by a derivative of the transition map along the transverse heteroclinic intersection of the cycle, which is denoted as $\alpha$ in  \citep[Equation (2.4)]{LiTur:24}.  See Section~\ref{sec:LT} for a brief explanation, and \citep[Section 2.3]{LiTur:24} for more details. Recall that two transitive hyperbolic sets are {\em homoclinically related} if they have the same index, and the unstable manifold of each set intersects the stable manifold of the other set transversely.
To facilitate the presentation, we call a transitive hyperbolic set  {\em $C^1$-wild} if  it exhibits uncountably many $C^1$-robust homoclinic tangencies, and use $\homo$ to denote the homoclinic relation.

 \begin{thm}\label{thm:hdc}
Let $f\in \diff^r(\mathcal{M})$, $r=1,\dots,\infty,\omega$,  have a coindex-1 heterodimensional cycle involving two saddles $O_1$ and $O_2$. Then, $f$ can be approximated in the $C^r$ topology by a diffeomorphism $g$ which has at least one standard blender, being simultaneously separated and arrayed. More specifically,
\begin{itemize}[nosep]
\item in the saddle case, 
\begin{itemize}[nosep]
\item  when $|\alpha|\leq 1$,  there is one such blender $\Lambda_1 \homo O_{1,g}$,
\item  when $|\alpha|\geq 1$, there is one such blender $\Lambda_2 \homo O_{2,g}$, and
\end{itemize}
\item  in other cases, there are two such blenders $\Lambda_1$ and $\Lambda_2 $ with $\Lambda_1\homo O_{1,g}$ and $\Lambda_2 \homo O_{2,g}$, where
\begin{itemize}
\item in saddle-focus  case,
\begin{itemize}
\item $\Lambda_1$ is $C^1$-wild if  $\gamma$ is nonreal,
\item $\Lambda_2$ is $C^1$-wild  if  $\lambda$ is nonreal, and
\end{itemize}
\item in the  double-focus case, $g$ can be chosen such that either $\Lambda_1$ or  $\Lambda_2$  is $C^1$-wild.
\end{itemize}
\end{itemize}
Moreover, for $i=1,2$, if we additionally assume that $O_i$ has a homoclinic tangency, then the blender $\Lambda_i$ is  always $C^1$-wild.
In all cases, $\Lambda_1$ is center-stable and  $\Lambda_2$ is center-unstable.
\end{thm}

Note that, in the saddle-focus case, both $\Lambda_1$ and $\Lambda_2$ are $C^1$-wild if either  $O_1$ has a homoclinic tangency and $\lambda$ is nonreal, or $O_2$ has a homoclinic tangency and $\gamma$ is nonreal; in the double-focus case, this happens if one of $O_1$ and $O_2$ has a homoclinic tangency.

Above two theorems follow from a more detailed result, Theorem~\ref{thm:hdcfull} in Section \ref{sec:LT}, which considers the coexistence of robust homoclinic tangencies and robust heterodimensional dynamics. Its proof is based on the analysis of the infinitely many standard blenders that coexist near a coindex-1 heterodimensional cycle \citep{LiTur:24} (i.e., there is a sequence of pairwise disjoint open sets near the cycle such that the locally maximal invariant subset of each open set is a standard blender). These coexisting blenders, all of which have nearly-affine central dynamics, are a generalization of the nearly-affine blender model introduced by \citep{BerCroPuj:22}. They form the {\em nearly-affine blender system} (see Definition \ref{defi:NABS}), which serves as a source for generating the blenders studied in this paper.


\subsection{Strengthening the Newhouse theorem}\label{sec:introNewhouse}

Theorem~\ref{thm:main}  strengthens  the classical Newhouse theorem (and its higher dimensional generalization) in the sense that the previous $C^2$-robust homoclinic tangencies are now replaced with $C^1$-robust ones, provided that the original tangency is associated with a saddle of effective dimension larger than two. The strategy of the proof is   to first obtain a heterodimensional cycle, by the result of \citep{LiLiShiTur:22}, and then apply Theorem~\ref{thm:hdc}.


Let us denote the multipliers of a saddle $O$
as 
\begin{equation*}\label{eq:eigenvalue}
|\lambda_{d_{s}}|\leqslant \dots\leqslant|\lambda_1|<1<|\gamma_1|\leqslant\dots\leqslant |\gamma_{d_u}|.
\end{equation*}
In general position, there is only one real or one pair of conjugate complex center-stable multipliers, namely, $\lambda_1$ or $\lambda_1$ and $\lambda_2=\lambda^*_1$. The same for center-unstable multipliers. Let us denote by $d_{cs}$ and $d_{cu}$ the numbers of center-stable and -unstable multipliers.
Suppose that $O$ has a generic homoclinic tangency (see e.g. \citep{LiLiShiTur:22} for the genericity conditions) and denote by $\Gamma$ the corresponding homoclinic orbit.
According to \citep{GonTurShi:93b,BonCro:16}, the orbit $\mathcal{O}(O)$ of $O$ and $\Gamma$  lie in a center invariant manifold $\mathcal{M}^c$ of dimension $(d_{cs}+d_{cu})$, which is tangent to the eigenspace corresponding to the central multipliers. 


\begin{defi}[Effective dimension]\label{defi:d_eff}
We say that a saddle has {\em effective dimension larger than one}  if $\dim \mathcal{M}^c >2 $  and the dynamics restricted to it is neither area-contracting nor area-expanding. More specifically, $d_{\rm eff}>1$ if
\begin{enumerate}[nosep]
\item $\lambda_1=\lambda^*_2$ is not real, $\gamma_1$ is real, and $|\lambda_1\gamma_1|>1$, or
\item $\lambda_1$ is real, $\gamma_1=\gamma^*_2$ is not real, and $|\lambda_1\gamma_1|<1$, or
\item 
 $\lambda_1=\lambda^*_2$  and $\gamma_1=\gamma^*_2$ are not real, and $|\lambda_1\gamma_1|\neq 1$.
\end{enumerate}
\end{defi}

 See \citep{Tur:96} for the precise definition of effective dimension (instead of ``effective dimension larger than one'').  It was established in  \citep{GonShiTur:08} that a homoclinic bifurcation associated with a saddle of effective dimension $d_\mathrm{eff}$ can give rise to   $d_\mathrm{eff}$ coexisting saddles of $f|_{\mathcal{M}^c}$ (and hence of $f$) with distinct indices. Intuitively, for instance in case~1 of Definition~\ref{defi:d_eff}, the condition on the multipliers implies  area expansion for $f|_{\mathcal{M}^c}$. This creates room for the birth of saddles  with index two, in addition to the original saddle $O$, which has index one  as a saddle of $f|_{\mathcal{M}^c}$. 
We therefore say that  the central dynamics -- namely, the dynamics restricted to $\mathcal{M}^c$ -- is {\em not essentially  two-dimensional}, in the sense that it allows the coexistence of saddles with different indices, which is impossible in dimension two.
The presence of these two saddles gives the possibility of creating a heterodimensional cycle involving them. This is indeed the case by  the following 
\begin{thm*}[{\citep[Theorem 1]{LiLiShiTur:22}}]
Let $f\in \diff^r(\mathcal{M})$, $r=1,\dots,\infty,\omega$,  have a homoclinic tangency to a saddle $O$ with $d_{\rm eff}>1$. Then, $f$ can be approximated in the $C^r$ topology by a diffeomorphism which has a coindex-1 heterodimensional cycle involving the continuation of $O$ and another saddle $Q$ satisfying   $\ind(Q)=\ind(O)+1$ if $|\lambda_1\gamma_1|>1$, and  $\ind(Q)=\ind(O)-1$ if $|\lambda_1\gamma_1|<1$. 
\end{thm*}

Here the assumption $|\lambda_1\gamma_1|>1$ means that $O$ belongs either to case 1 of Definition \ref{defi:d_eff}, or to case 3  with $|\lambda_1\gamma_1|>1$. Similarly for the assumption $|\lambda_1\gamma_1|<1$.

By the (multidimensional) Newhouse theorem, up to an arbitrarily $C^r$-small perturbation, one can assume that $f$ has a hyperbolic basic set having a $C^2$-robust homoclinic tangency. Moreover, this set is homoclinically related to $O$  \citep{GonTurShi:93c,PalVia:94}. Hence, after applying the perturbation involved in the above theorem, one can, by the lambda lemma, recover the homoclinic tangency to the continuation of $O$,  by an additional arbitrarily $C^r$-small perturbation. After that, invoking the saddle-focus case of Theorem \ref{thm:hdc} immediately gives
\begin{thm}\label{thm:ht}
Let $f\in \diff^r(\mathcal{M})$, $r=1,\dots,\infty,\omega$, have a homoclinic tangency to a saddle  $O$ with $d_{\rm eff}>1$. Then, $f$ can be approximated  in the $C^r$ topology by a diffeomorphism $g$ which has a standard cs-blender and a standard cu-blender, $\Lambda_1$ and $\Lambda_2$ with $\ind(\Lambda_1)+1=\ind(\Lambda_2)$, each of which is simultaneously separated and arrayed, and exhibits uncountably many $C^1$-robust homoclinic tangencies. Moreover,   $O_g\homo \Lambda_1$  if $|\lambda_1\gamma_1|>1$, and   $O_g\homo\Lambda_2$ if $|\lambda_1\gamma_1|<1$. 
\end{thm}

\begin{rem}\label{rem:nece}
The condition  $d_\eff>1$ is   {\em necessary} for creating heterodimensional cycles  from the bifurcation of a single\footnote{Heterodimensional cycles can also be created in case 1 of Definition \ref{defi:d_eff} if $O$ has two orbits of homoclinic tangency which  intersect the same strong-stable/unstable leaf, see \citep{LiTur:20}.}
generic homoclinic tangency, and, hence, necessary for creating $C^1$-robust homoclinic tangencies through the mechanism of heterodimensional bifurcations. To see this, note that any heterodimensional cycle born from the  homoclinic bifurcation of $O$ is also a cycle for $f|_{\mathcal{M}^c}$. So, the involved hyperbolic periodic points, when viewed as those of $f|_{\mathcal{M}^c}$, must also  be saddles of different indices, which is impossible when $d_\eff= 1$, see~\citep{LiLiShiTur:22} for details.
\end{rem}

The robust presence of uncountably many homoclinic tangencies also appears in the two-dimensional Newhouse construction (though only $C^2$-robust). This was not proved in his original paper, but follows from the main result in \citep{HunKanYor:93}. For a detailed explanation see \citep[Corollary 2.17]{BerBie:22}.

It was shown by Gonchenko \citep{Gon:83}, generalizing the two-dimensional result by Gavrilov-Shilnikov \citep{GarShi:73} and Newhouse \citep{New:74}, that a sectionally-dissipative saddle $O$ with a homoclinic tangency can produce sinks by an arbitrarily small perturbation. Thus, in any open set where diffeomorphisms with a homoclinic tangency to a sectionally-dissipative saddle are dense, 
there is a generic subset where every map displays infinitely many sinks \citep{New:79, Gon:83,PalVia:94}. Later, Gonchenko, Shilnikov and Turaev \citep{Tur:96,GonShiTur:08} relaxed the sink-producing condition   such that, instead of sectional dissipation, it  requires the absence of strong expansion and the contraction of volumes in the center manifold $\mathcal{M}^c$. 
\begin{defi}[Saddles of sink-producing type]\label{defi:sink}
The saddle $O$  is said to be  of  {\em sink-producing type} if 
\begin{enumerate}[nosep]
\item  $\lambda_1=\lambda^*_2$ is not real, $\gamma_1$ is real, $|\lambda_1^2\gamma_1|<1$, and $\dim(W^u(O))=1$,
\item  $\lambda_1$ is real, $\gamma_1=\gamma^*_2$ is not real,  $|\lambda_1\gamma_1^2|<1$, and $\dim(W^u(O))=2$,
\item 
 both $\lambda_1=\lambda^*_2$  and $\gamma_1=\gamma^*_2$ are not real, $|\lambda_1\gamma_1|< 1$, and $\dim(W^u(O))=2$, or
\item  both $\lambda_1$ and $\gamma_1$ are real,  $|\lambda_1\gamma_1|<1$ and $\dim(W^u(O))=1$.
\end{enumerate}
\end{defi}

\begin{proof}[Proof of Corollary~\ref{cor:main}]
By Theorem~\ref{thm:ht}, $f$ can be $C^r$-approximated by a diffeomorphism $g$ which has a hyperbolic basic set $\Lambda$ exhibiting a $C^1$-robust homoclinic tangency. Let $\mathcal{U}$ be the neighborhood of $g$ associated with the robust tangency.   Since $\Lambda$ is homoclinically related to the continuations of $O$, it follows from the lambda-lemma that there exists a $C^r$-dense subset $\mathcal{U}'$ of $\mathcal{U}$ such that, for every diffeomorphism $h\in \mathcal{U}$, the continuation $O_h$ has a homoclinic tangency. Moreover, up to shrinking $\mathcal{U}$, the saddle $O_h$ is of sink-producing type if $O$ is. Since \citep[Theorem 5]{GonShiTur:08} asserts that sinks can be born from sink-producing saddles by an arbitrarily $C^r$-small perturbation, one  easily obtains the corollary from a standard genericity argument (see e.g. \citep{New:79}).
\end{proof}

\subsection{Mechanism for generating robust tangencies  and comparison of different types of blenders}\label{sec:differ}
 We begin by briefly describing the structure of a blender-horseshoe. To simplify the presentation, we restrict to dimension three.
A  blender-horseshoe $\Lambda$ of a diffeomorphism $f$ is characterized by three defining features. Firstly, it is a hyperbolic basic set   with a partially hyperbolic structure\footnote{The blender-horseshoe introduced in \citep{BonDia:12b} has an expanding center. Here, in order to be consistent with the constructions in the present paper, we consider instead a contracting center. The two settings are equivalent up to replacing $f$ by $f^{-1}$.}, that is, the tangent bundle restricted to $\Lambda$ admits an invariant splitting $T \Lambda=E^{ss}\oplus E^{cs}\oplus E^{u}$; in the three-dimensional case, all subbundles are one-dimensional, and vectors in $(E^{ss}\oplus E^{cs})$ and $E^{u}$  are uniformly contracted and, respectively, expanded by the differential $Df$.  There are coordinates such that $\Lambda$ has a cube $\Pi \subset \mathbb{R}^3$ as its isolating neighborhood (i.e., $\Lambda=\bigcap_{n\in \mathbb{Z}}f^n(\Pi)$). Moreover, $E^{ss}$ lies in a $D(f^{-1})$-invariant  constant cone field  $\mathcal{C}^{ss}$ around the $z$-direction and $E^{u}$ lies in a $Df$-invariant  constant cone field  $\mathcal{C}^{u}$ around the $y$-direction. In particular, every local unstable leaf $\ell^u\subset W^u(\Lambda)\cap \Pi$ is tangent to $\mathcal{C}^{u}$.
These cones are sufficiently small so that, if a curve is tangent to $\mathcal{C}^{ss}$ or $\mathcal{C}^{u}$, then it is nearly parallel to the $z$-axis or $y$-axis and crosses the cube $\Pi$.  

The second defining feature of $\Lambda$, which distinguishes it from other hyperbolic sets, is the ``blender property'': for every curve $S$ that is tangent to $ \mathcal{C}^{ss}$ and crosses $\Pi$, one has $S\cap \ell^u \neq \emptyset$ for some unstable leaf $\ell^u$ (see Proposition~\ref{prop:nontr}). The last defining feature is the so-called ``reproducing property", which makes a blender-horseshoe special among general  blenders. In what follows, we describe this property and its role in the creation of robust tangencies, following the original idea in \citep{BonDia:12b}.

\begin{figure}[h]
\begin{center}
\includegraphics[scale=1.1]{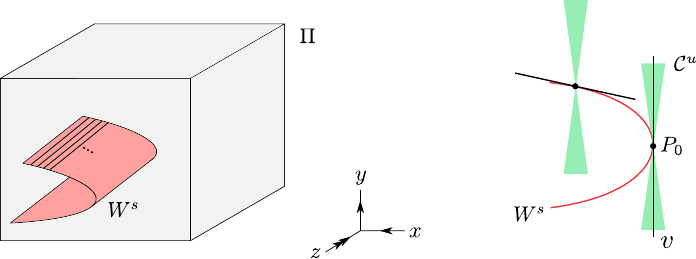}
\caption{Left picture: a parabola-like piece of $W^s(\Lambda)$ that crosses $\Pi$ and is foliated by curves tangent to $\mathcal{C}^{ss}$. Right picture: the tangent space $T_{P_0}W^s$ contains a vector $v\in \mathcal{C}^{u}$.}
\label{fig:mechanism}
\end{center}
\end{figure}

Now suppose that $\Lambda$ has a quadratic homoclinic tangency. One can  find a piece $W^s$ of $W^s(\Lambda)$ which intersects non-transversely some unstable leaf $\ell^u\subset \Pi$ and  is foliated by curves almost parallel to the $z$-axis, that is, $W^s=\bigcup S$ for a collection of curves $S$ that are tangent to $\mathcal{C}^{ss}$ and cross $\Pi$ (see Figure~\ref{fig:mechanism}).  By the blender property, every  such curve $S$ intersects  some   leaf $\ell^u$. Since all $\ell^u$ are tangent to $\mathcal{C}^u$, it follows from the quadraticity of the tangency that there exists a point $P_0$ in some unstable leaf such that $T_{P_0} W^s\cap \mathcal{C}^u_{P_0}\neq 0$ (see Figure~\ref{fig:mechanism}).  The blender-horseshoe $\Lambda$ has the ``reproducing property'': the preimage of $W^s$ by $f$ contains a piece which is also foliated by curves  that are tangent to $\mathcal{C}^{ss}$ and cross $\Pi$ (see Lemma~\ref{lem:folding2}). Iterating $W^s$ infinitely many times then yields a sequence $\{P_n\}$ of points such that $T_{P_n}f^{-n}(W^s)\cap \mathcal{C}^u_{P_n}\neq 0$. One then easily finds  $T_{P_*}W^s\cap T_{P_*} \ell^u_* \neq 0$, where $\ell^u_*$ is the unstable leaf containing the limit point $P_*=\lim P_n$. This gives the desired homoclinic tangency of $\Lambda$ (see Proposition~\ref{prop:folding}). Since the property of being a blender-horseshoe is  $C^1$-robust, the tangency is  $C^1$-robust as well.

The main constraint in the existing construction (e.g. \citep{BonDia:12b,BerCroPuj:22}) of blender-horseshoes is that the contraction/expansion  in the central direction is close to $1$. In this case, the action of $f^{-1}$ on $\Pi$ induces an iterated function system (IFS) of two elements, and it is relatively easy to control the preimages of $W^s$ and establish the ``reproducing property''. However, this constrain makes it difficult to obtain blender-horseshoes in a general setting -- it is one of the main reasons why blender-horseshoes are obtained in \citep{BonDia:12b} via $C^1$, rather than $C^r$, perturbations of heterodimensional cycles.

The standard blenders constructed in \citep{LiTur:24} remove the restriction on the center and can be obtained through $C^r$ perturbations near heterodimensional cycles; however,  only the ``blender property'' is guaranteed.  The separated blenders defined in the present paper are a special class of standard blenders: they not only allow arbitrary expansion or contraction rates in the central direction (so the IFS induced from $f^{-1}$ may contain  a large number of elements), but also retain the ``reproducing property'' that is crucial for the generation of $C^1$-robust tangencies. 
The arrayed blenders are even stronger, as they allow the ``reproduction'' of several copies simultaneously at each iteration and lead to the creation of uncountably many robust tangencies. 


We  note that the structure of separated  blenders have the potential to be extended to the case where the central dimension is larger than one, and hence be used to create robust  tangencies of a {\em higher corank}, that is, the tangent spaces of the stable and unstable manifolds of the blender at the tangency point share a common subspace of dimension larger than one, see Remark~\ref{rem:diff}. For recent progresses on high corank tangencies, see \citep{BarRai:17,BarRai:21,Asa:22,Min:26}.
 Although arrayed blenders   seem to be a more powerful version of the separated ones, the construction of arrayed blenders is fundamentally different, as it   relies on  the one-dimensional center in an essential way; hence they cannot be easily generalized to the case of higher-dimensional centers. 

It is worth mentioning that  arrayed  blenders can  be embedded into parameter families such that they become the parablenders introduced by Berger \citep{Ber:16}, which are key to build an open region, the {\em Berger domain}, in the space of finite-parameter families of diffeomorphisms where generic families display infinitely many sinks for an open set of parameter values. See \citep{BerCroPuj:22,BarRai:24} for  recent advances in this direction. 
In a forthcoming paper \citep{LiTur:xxb}, we use arrayed  blenders to show that  Berger domains can be found near every homoclinic tangency of the sink-producing type.\\

In the next three sections, we define standard blenders, separated standard blenders, and arrayed standard blenders, respectively, which are characterized by Markov partitions, and prove their fundamental properties on the creation of $C^1$-robust homoclinic tangencies.
In Section \ref{sec:hdc}, we prove Theorem~\ref{thm:hdc} by showing that separated and arrayed blenders exist near heterodimensional cycles.

\section{Standard blenders}\label{sec:blender}
The notion of blenders was introduced by Bonatti and D\'iaz for constructing robust non-hyperbolic transitive diffeomorphisms \citep{BonDia:96}. A blender is a hyperbolic  basic  set with the property that the projection of its stable or unstable lamination to certain central subspace has a non-empty interior. As a consequence, the lamination acts as if it has one more topological dimension and hence can produce robust non-transverse intersections with other submanifolds. 

Blenders have been used to obtain important dynamical phenomena in a wide range of settings. For a list of these applications  we refer the readers to \citep{Li:24} and the references therein. Due to different purposes, blenders have emerged in various forms during the last decades, see, e.g. \citep{NasPuj:12,BarKiRai:14,BocBonDia:16,AviCroWil:21}. In particular, the blenders arisen from heterodimensional cycles are characterized by Markov partitions of a finite (and possibly large) number of elements, which motivated the notion of  standard blenders given in \citep[Appendix]{LiTur:24}. In comparison, a blender-horseshoe in \citep{BonDia:12b} corresponds to a Markov partition of exactly two elements. We note that the idea of Markov partitions is also used (implicitly) in \citep{CapKraOsiZgl:23} to achieve a characterization of blenders which is in particular  suitable for rigorous computer-assisted validation. 

In what follows, we first introduce partially hyperbolic Markov partitions  and define standard blenders. After that, we prove    the fundamental property of standard blenders  that they can produce robust non-transverse heterodimensional intersections, see Proposition \ref{prop:nontr}.

%
%

\subsection{Partially hyperbolic Markov partitions}\label{sec:markov}

Let $\{\Pi_i\}$ be a finite collection of pairwise disjoint closed subsets of a manifold  $\mathcal{M}$ of dimension $d\geq 3$, and $g:  U\to \mathcal{M} $ be any diffeomorphism of some neighborhood $U$ of $\bigcup  \Pi_i$ such that $g(\Pi_i)\subset U$.
Up to dividing $\Pi_i$ into subsets, we can assume that the intersection $\Pi_i\cap g(\Pi_j)$ has at most one connected component.  We call the  pair $(\{\Pi_i\},g)$ a {\em Markov partition} if  for each pair $(i,j)$ it holds that $g^n(\Pi_i)\cap \Pi_j\neq \emptyset$ for some $n\geq 0$. 
Denote by $\inte(\cdot)$ the interior of a set.

\begin{defi}[cs-Markov partitions]\label{defi:markov}
A   Markov partition $(\{\Pi_i\},g)$ is said to be {\em center-stable (cs) partially hyperbolic} if
\begin{itemize}[nosep]
\item (Horizontal domains) there exist charts $(U,\pp_i)$, convex closed subsets $\mathbb{X}_i\subset \mathbb{R}, \mathbb{Y}_i\subset \mathbb{R}^{d_u},\mathbb{Z}_i\subset \mathbb{R}^{d_{ss}}$ with $1+d_u+d_{ss}=d$,  and smooth functions $\psi_{i}:\mathbb X_i\times \mathbb R^{d_u} \times \mathbb Z_i \to \inte(\mathbb Y_i)  $ such that
\begin{equation}\label{eq:phipi}
\pp_i(\Pi_i)=\{(x,\psi_{i}(x,y,z),z):(x,y,z)\in \mathbb X_i\times \mathbb Y_i \times \mathbb Z_i \};
\end{equation}
\item (Hyperbolicity) for every pair $(i,j)$ with $g(\Pi_i)\cap\Pi_j\neq \emptyset$, there exists a function $g^\times_{ij}:\mathbb X_i\times \mathbb Y_j \times \mathbb Z_i \to \inte(\mathbb X_j\times \mathbb Y_i \times \mathbb Z_j)$ such that 
 $\pp_j\circ g\circ \pp^{-1}_i(x,y,z)=(\bar x,\bar y,\bar z)\in \mathbb X_j\times \mathbb Y_j \times \mathbb Z_j$ if and only if 
 \begin{equation}\label{eq:cross}
(\bar x,y,\bar z)=g^\times_{ij}( x, \bar y,z),
\end{equation}
and, with denoting  $g^\times_{ij}=(g^\times_{ij1},g^\times_{ij2},g^\times_{ij3})$,   it holds for some $ \mu \in (0,1)$ that 
\begin{equation}\label{eq:deri1}
\left\| \dfrac{\partial (g^\times_{ij1},g^\times_{ij3})}{\partial (x,z)}\right\|
+\left\| \dfrac{\partial (g^\times_{ij1},g^\times_{ij3})}{\partial \bar y}\right\|<\mu,\qquad
\left\| \dfrac{\partial g^\times_{ij2}}{\partial (x,z)}\right\|
+\left\| \dfrac{\partial g^\times_{ij2}}{\partial \bar y}\right\|<\mu,
\end{equation}
in some suitable norm;
\item (Partial hyperbolicity) one has  
$$\left| \frac{\partial g^\times_{ij1}}{\partial x }\right|\in (0,1),$$ 
 and the function  $\tilde g^\times_{ij}:\mathbb X_j\times \mathbb Y_j \times \mathbb Z_i \to \inte(\mathbb X_i\times \mathbb Y_i \times \mathbb Z_j)$, obtained from
 $ g^\times_{ij}$ by inverting\footnote{
First express $x$ as a function of $(\bar x,\bar y,z)$ through $\partial g^\times_{ij1}$, and then substitute this expression into $g^\times_{ij2}$ and $g^\times_{ij3}$.
}
 its first coordinate  $g^\times_{ij1}$, satisfies, for some $\nu\in (0,1)$, that
\begin{equation}\label{eq:deri2}
\begin{gathered}
\left\| \dfrac{\partial \tilde g^\times_{ij3}}{\partial z}\right\| +
\left\| \dfrac{\partial \tilde g^\times_{ij3}}{\partial (\bar x,\bar y)}\right\|
<1,\\
\left\| \dfrac{\partial (\tilde g^\times_{ij1},\tilde g^\times_{ij2})}{\partial (\bar x,\bar y)}\right\|
\cdot
\left\| \dfrac{\partial \tilde g^\times_{ij3}}{\partial z}\right\|
\cdot
\left(1- \left\| \dfrac{\partial \tilde g^\times_{ij3}}{\partial (\bar x,\bar y)}\right\|\right)^{-1}
+
\left\| \dfrac{\partial (\tilde g^\times_{ij1},\tilde g^\times_{ij2})}{\partial z}\right\|
<\nu.
\end{gathered}
\end{equation}
\end{itemize}
\end{defi}

\begin{rem}
Similarly defined Markov partitions are also considered in \citep{Tur:24} in a hyperbolic setting,  without assuming partial hyperbolicity.
\end{rem}

Let us make some simple observations from the above definition.
To simplify the notations, we will always identify $\Pi_i$ with $\pp_i(\Pi_i)$ throughout Sections \ref{sec:blender}--\ref{sec:array}, unless otherwise stated.

By \eqref{eq:phipi} and \eqref{eq:cross}, each preimage $H_{ij}:=\Pi_i\cap g^{-1}(\Pi_j)$ is a  {\em horizontal strips} in the sense that it has full $(x,z)$-size and small $y$-size in $\Pi_{i}$. That is, $H_{ij}$ has a skew-product structure as in \eqref{eq:phipi}, but with some function $\psi_{ij}$ satisfying $|\psi_{ij}(x,y,z)|<|\psi_{i}(x,y,z)|$. Each image $V_{ji}:=\Pi_i\cap g(\Pi_{j})$ is a {\em vertical strip} in the sense that $g(\Pi_j)$ has full $y$-size and small $(x,z)$-size in $\mathbb X_i\times \mathbb Y_i \times \mathbb Z_i$ by \eqref{eq:cross}, and hence  $V_{ji}$ has full $y$-size and small $(x,z)$-size in $\Pi_i$ by \eqref{eq:phipi}.  

Conditions \eqref{eq:deri1} and \eqref{eq:deri2} mean that the $x,y$, and $z$ directions are mapped approximately  to the $\bar x,\bar y$, and $\bar z$ directions, respectively; moreover, there is weak contraction in $x$, strong contraction in $z$, and expansion in $y$. More precisely, using \eqref{eq:deri1} and \eqref{eq:deri2}, one can show\footnote{For an illustration of such cone computations, see e.g. \citep[Lemma~3.1]{LiTur:24}
} 
 the existence of the following cone fields on  $\Pi_i$ (where $(dx,dy,dz)$ denotes a tangent  vector):
\begin{equation}\label{eq:cone}
\begin{aligned}
\mathcal{C}^{u}&=\{(dx,dy,dz):|dx|+\|dz\|\leq K_i^u |dy|\},\\
\mathcal{C}^{s}&=\{(dx,dy,dz):\|dy\|\leq K_i^s (|dx|+\|dz\|)\},\\
\mathcal{C}^{ss}&=\{(dx,dy,dz):|dx|+\|dy\|\leq K_i^{ss} \|dz\|\},
\end{aligned}
\end{equation}
for some $K_i^u,K_i^s \in (0,1)$, depending on $\mu$, and $K_i^{ss} \in (0,1)$, depending on $\nu$, such that, for any $P\in \Pi_i $ and $\bar P=g(P) \in \Pi_j$,
\begin{itemize}[nosep]
\item $\mathcal{C}^u$ is forward-invariant:  $\mathcal{C}^u_P$ is mapped strictly inside $\mathcal{C}^u_{\bar P}$ by the differential $D(\pp_j\circ g\circ \pp^{-1}_i)$,
\item $\mathcal{C}^s$ and $\mathcal{C}^{ss}$ are backward-invariant:  $\mathcal{C}^{s/ss}_{\bar P}$ is mapped strictly inside $\mathcal{C}^{s/ss}_{P}$ by $D(\pp_i\circ g^{-1}\circ \pp^{-1}_j)$, 
\item vectors in $\mathcal{C}^{u}$ are uniformly expanded by $D(\pp_j\circ g\circ \pp^{-1}_i)$ and vectors in $\mathcal{C}^{s/ss}$ are uniformly contracted by $D(\pp_i\circ g^{-1}\circ \pp^{-1}_j)$.
\end{itemize}


\begin{ex}\label{ex}
A simple example of a cs-Markov partition can be generated by a three-dimensional horseshoe, as illustrated in  Figure \ref{fig:3dhorseshoe}. The map $g$ contracts the cube in the $(x,z)$-directions, expand it in the $y$-direction, and then bend it such that its top crosses the cube along the $y$-direction. This results in a cs-Markov partition $\{\Pi_1,\Pi_2\}$, and the above-mentioned horizontal and vertical strips are in fact  substrips inside $\Pi_i$ and $g(\Pi_i)$.
\begin{figure}[!h]
\begin{center}
\includegraphics[scale=.9]{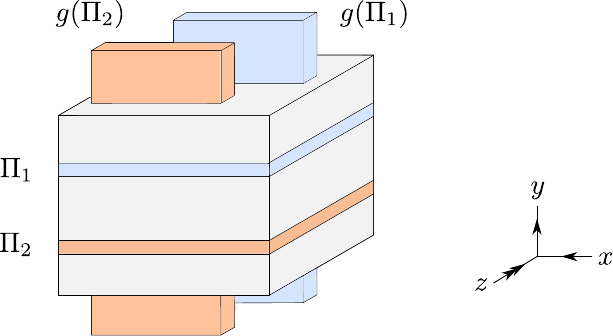}
\caption{A three-dimensional Horseshoe}
\label{fig:3dhorseshoe}
\end{center}
\end{figure}
\end{ex}

%
%
%

\begin{conv}
For a Markov partition $(\{\Pi_i\}_{i=1}^N,g)$, iterations of $g$ are defined inductively by $g^{n+1}(\Pi_i)=\bigcup_{j=1}^N g(g^n(\Pi_i)\cap\Pi_j)$.
\end{conv}

Denote $\Pi=\bigcup \inte(\Pi_i)$. The Markov partition induces a zero-dimensional hyperbolic basic set $\Lambda:=\bigcap_{n\in\mathbb Z} g^{n}(\Pi)$ consisting of points whose orbits never leave $\Pi$.  To see this, note first that, since hyperbolicity follows directly from the existence of the cone fields, it suffices to show that this invariant set is non-empty. By definition, the map $g^\times_{ij}$ takes $\mathbb X_i\times \mathbb Y_j\times \mathbb Z_i$ into $\inte(\mathbb X_j\times \mathbb Y_i\times \mathbb Z_j)$ and is a  contraction. It then follows from a fixed-point lemma on products of  metric spaces (see e.g. \citep[Theorem 5.8]{Shi:67}) that the set  of points whose orbits by $g$ lie entirely in $\Pi$, i.e., the set $\Lambda$, is in one-to-one correspondence to the set $\{1,\dots,\card \{\Pi_i\}\}^{\mathbb Z}$, or a subset of it when there exist $i$ and $j$ such that $g(\Pi_i)\cap \Pi_j=\emptyset$. More specifically, each point $P\in\Lambda$ corresponds to a coding $(i_n)_n\in \{1,\dots,\card \{\Pi_i\}\}^{\mathbb Z}$ such that its orbit $\{P_n:=g^n(P)\}$  satisfies $P_n\in \Pi_{i_n}$.  In particular,  $\Lambda$  is   non-empty. By construction, the stable manifold $W^s(\Lambda)$ consists of {points whose forward iterates never leave $\Pi$}, while the unstable manifold $W^u(\Lambda)$ consists of {points whose backward iterates never leave $\Pi$}.



\begin{defi}[cu-Markov partitions]\label{defi:Markovcu}
We say that the  $(\{\Pi_i\},g)$ is a  {\em center-unstable partially hyperbolic (cu)}  Markov partition if its  inverse  $(\{g(\Pi_i)\},g^{-1})$ is a cs-Markov partition.
\end{defi}

\subsection{Robust heterodimensional non-transverse intersections}\label{sec:hetero}
The complexity of the dynamics arisen from a cs-Markov partition, and, in particular, our results, are tightly related to the behavior of the backward iteration of curves tangent to the cone field $\mathcal{C}^{ss}$, which we call {\em ss-discs} and define below. We are interested the recurrence property of such discs that the preimage of every ss-disc contains a new ss-disc.

Let  $\partial_x \Pi_i$ be  the $x$-boundary of $\Pi_i$, namely, the set  $ \partial \mathbb X_i \times \mathbb Y_i \times\mathbb Z_i$. We also have the $y$-boundary $\partial_y \Pi_i$ and $z$-boundary $\partial_z \Pi_i$. Similarly, we define the $x,y,z$-boundaries of horizontal strips $H$ of $\Pi_i$. For a vertical strip of the form $V=g(H)\cap \Pi_i$ for some horizontal strips $H\subset \Pi_j$, we have $\partial_\rho V=g(\partial_\rho H)\cap \Pi_i$ for $\rho=x,y,z$.

\begin{defi}[ss-discs and u-discs]\label{defi:disc}
Let $\mathcal{C}^{ss}$ and $\mathcal{C}^{u}$ be the cone fields in \eqref{eq:cone} associated to the cs-Markov partition $(\{\Pi_i\},g)$.
An {\em ss-disc} $S$ of $\Pi_i$ is a submanifold tangent to the cone field $\mathcal{C}^{ss}$ such that it lies in $ \Pi_i\setminus (\partial_x \Pi_i\cup \partial_y \Pi_i)$ and is given by some smooth function  $(x,y)=\s(z)$ defined for all $z\in\mathbb Z_i$. 
A {\em u-disc} $S^u$ of $\Pi_i$ is a submanifold tangent to the cone field $\mathcal{C}^{u}$ such that it lies in $ \Pi_i\setminus (\partial_x \Pi_i\cup \partial_z \Pi_i)$ and is  given by some smooth function  $(x,z)=\s^u(y)$ defined for all $y\in\mathbb Y_i$. 
\end{defi}


\begin{lem}\label{lem:cover}
Let $S$ be an ss-disc of $\Pi_j$ and $V$ be a vertical strip of the form $V=g(H)$ for some horizontal strip $H\subset \Pi_i$.
If $S$ crosses $V$, i.e., $S\cap V\subset (V\setminus \partial_x V)$, then $g^{-1}(S\cap V)$ is an ss-disc of $\Pi_i$.
\end{lem}


\begin{proof}
The intersection $V\cap \pi_i$ is given by the solution to the system consisting of the function $\tilde g^\times_{ij}$ in Definition~\ref{defi:markov}:
$$x= \tilde g^\times_{ij1}(\bar x,\bar y,z),\qquad
y= \tilde g^\times_{ij2}(\bar x,\bar y,z),\qquad
\bar z= \tilde g^\times_{ij3}(\bar x,\bar y,z), \
$$
defined for $(\bar x,\bar y,z) \in \X_j \times \Y_j \times \Z_i $, and the defining function  of $S$: 
$$\bar x =s_1(\bar z),\qquad \bar y=s_2(\bar z), $$
defined for $\bar z\in \Z_j $ and satisfying $\|\partial s_{1,2}/\partial \bar z\| <1 $ by \eqref{eq:cone}. Note that the first line of~\eqref{eq:deri2} implies $\|\partial \tilde g^\times_{ij3}/\partial (\bar x,\bar y)\|<1$. Hence, by substituting $\bar z= \tilde g^\times_{ij3}(\bar x,\bar y,z)$ into the above defining function, we can find $(\bar x,\bar y)=\tilde s(z)$ for some function $\tilde s$ defined near a point in $S\cap V$. Since $S\cap V\subset (V\setminus \partial_x V)$, the function $\tilde s$ can be extend to $\Z_i$. Combining $\tilde s$ with  $\tilde g^\times_{ij1}$ and $\tilde g^\times_{ij2}$, we find a disc crossing $\Pi_i$. It is further an ss-disc of $\Pi_j$ due to the 
invariance of the cone field $\mathcal{C}^{ss}$.
\end{proof}

\begin{defi}[Base of $\Pi_i$]\label{defi:base}
We associate to each $i$ an open interval $\mathbb{X}_i^\mathrm{B}\subset \mathbb{X}_i$ and define
\begin{itemize}[nosep]
\item {\em base}: $\Pi^{\mathrm{B}}_i
=\Pi_i\cap(
{\mathbb{X}}^{\mathrm{B}}_i\times\mathbb Y_i\times\mathbb Z_i), 
\quad H^{\mathrm{B}}=H\cap \Pi^{\mathrm{B}}_i,
\quad V^{\mathrm{B}}=g(\tilde H^{\mathrm{B}})\cap \Pi^{\mathrm{B}}_i$,
\end{itemize}
where $H$ and $V$ are any horizontal and vertical strips of $ \Pi_i$, with $V=g(\tilde  H)\cap \Pi_i$ for some $\tilde H$ of $\Pi_j$.
We say that an ss-disc $S$ {\em crosses} the base $\Pi_i^{\mathrm{B}}$ if $S\cap \Pi_i\subset \Pi_i^{\mathrm{B}}$; similarly    for all horizontal strips $H$ and vertical strips $V$.
\end{defi}

\begin{defi}[Covering property]\label{defi:cover}
A cs-Markov partition $(\{\Pi_i\},g)$
has  the {\em  covering property} if
\begin{itemize}[nosep]
\item[\setword{(A1)}{word:A1}] for all $\Pi_i$, every ss-disc $S$ crossing the base $\Pi^\mathrm{B}_i$  crosses the base $V^\mathrm{B}$ of some vertical strip. 
\end{itemize} 
\end{defi}

This property is easy to achieve. Since $V$ has full $y$-size, whether $S$ intersects $V$ essentially depends on the $x$-coordinates of its points. When  the cone constants in  \eqref{eq:cone} are small, one can find, for some reference point $z^*_i$, an interval $I$ such that an ss-disc crosses $V$ whenever $\s_x(z^*)\in I$, where $\s_x$ is the $x$-component of its defining function. The covering property is essentially the existence of sufficiently many vertical strips whose corresponding intervals $I$ cover $\mathbb X^\mathrm{B}_i$. 

\begin{ex}\label{ex:base}
Continue with Example \ref{ex}, we define the bases by  darker colors as illustrated in Figure \ref{fig:base}. The covering property is satisfied since every ss-disc crossing $\Pi^\B_i$ must cross the one of the two intersections $\Pi^\B_i\cap g(\Pi^\B_1)$  and $\Pi^\B_i\cap g(\Pi^\B_2)$, each of which is the base of a vertical strip.
\begin{figure}[!h]
\begin{center}
\includegraphics[scale=.9]{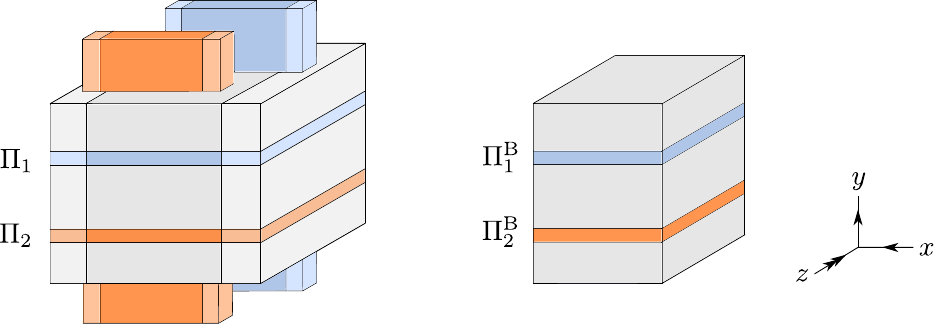}
\caption{$\Pi_i$ with bases}
\label{fig:base}
\end{center}
\end{figure}
\end{ex}

\begin{defi}[Standard blenders]\label{defi:blender}
Given a cs-Markov partition $(\{\Pi_i\},g)$ that satisfies the covering property, the locally maximal invariant set $\Lambda:=\bigcap_{n\in\mathbb Z} g^{n}(\bigcup \inte(\Pi_i))$ is called  a {\em standard center-stable (cs)  blender}. It is called  a {\em standard center-unstable (cu)  blender} if the Markov partition is center-unstable.
\end{defi}


As discussed after Definition~\ref{defi:markov}, a standard blender is a zero-dimensional  hyperbolic basic set, restricted to which $g$ is  conjugate to a full shift of $\card \{\Pi_i\}$ symbols. 
For cu-Markov partitions, one has, instead of~\eqref{eq:cone}, forward-invariant cone fields $\mathcal{C}^{u},\mathcal{C}^{uu}$ and backward-invariant cone field $\mathcal{C}^{s}$, which lead to the notions of  uu- and s-discs, completely parallel to the ss- and u-discs in Definition~\ref{defi:disc}.

\begin{prop}[Robust non-transverse intersections]\label{prop:nontr}
The unstable manifold of a standard cs-blender intersects every ss-disc which crosses some $\Pi^\B_i$; and the stable manifold of a standard cu-blender intersects every uu-disc  which crosses some $\Pi^\B_i$. 
\end{prop}
This is the fundamental property of blenders as introduced in the discovery paper \citep{BonDia:96} of blenders.
\begin{proof}
We only prove the proposition for cs-blenders, as the arguments for cu-blenders are parallel. Let $\Lambda$ be the blender induced from a Markov partition $(\{\Pi_i\},g)$. By the covering property, any ss-disc $S$ crossing $\Pi^\mathrm{B}_i$ crosses some $V^\mathrm{B}$. In particular, it crosses $V$, and hence, by  Lemma~\ref{lem:cover}, the preimage $S_1:=g^{-1}(S\cap V)$ is an ss-disc of some $\Pi_{i_1}$. By construction, $S_1$ crosses $\Pi^\B_{i_1}$. Repeating this process yields a coding $\underline{i}=(i_n)$ and a sequence of ss-discs $S_n$ crossing $\Pi^\B_{i_n}$. The images $\hat S_n:=g^{n}(S_n)$ are nested compact sets in $S$, so one can find a point $P\in \bigcap \hat S_n$ such that its backward iterates never leave $\bigcup \Pi_i$. Thus, $P\in W^u(\Lambda)$  by definition of $W^u(\Lambda)$, and, in particular, it belongs to the unstable leaf $\ell^u:=\bigcup_{n\geq 0} g^n(\Pi^\B_{i_n})$.  (See the discussion in the end of Section~\ref{sec:markov}.)
\end{proof}

\begin{rem}\label{rem:unstableleaf}
It follows from the proof that if an ss-disc crosses some $V^\mathrm{B}$, then it must intersect an unstable leaf that belongs to $V^\B$.
\end{rem}


\section{Separated standard blenders}\label{sec:sepa}
It is established  in \citep[Theorem D]{LiTur:24} (and also in \citep[Corollary A]{Li:24} as an improved version)  that a heterodimensional cycle produces an infinite sequence of standard blenders.  
An important observation is that an appropriate superposition of finitely many standard blenders from this system can yield more powerful ones, for example, the separated and arrayed blenders. 

As mentioned, these two types of blenders both possess the “reproducing property,” which is crucial for generating robust tangencies (see Section \ref{sec:differ}), but their constructions are fundamentally different.
In simple terms, separated blenders rely on the ability to distinguish the ``inside'' and ``outside'' of $\Pi_i$  in the central direction, whereas arrayed blenders require a notion of ``left'' and ``right'' of $\Pi_i$. While arrayed blenders are stronger in the sense of producing uncountably many robust tangencies (and give rise to Berger parablenders), they are restricted to the case of a one-dimensional center. By contrast, the notion of ``inside'' and ``outside'' can be generalized to higher dimensions, making it possible to construct  separated blenders with  center dimension higher than one, which can be used to obtain robust heterodimensional cycles of higher indices and robust homoclinic tangencies of  higher corank  (see Remark~\ref{rem:diff} for a further discussion).

\subsection{Base-Center structure}
We introduce an additional structure for a cs-Markov partition $(\{\Pi_i\},g)$. Recall that we have identified $\Pi_i$ with $\mathbb X_i\times\mathbb Y_i\times\mathbb Z_i$.

\begin{defi}[Base-Center structure of $\Pi_i$]\label{defi:structure1}
 One can associate to each $i$ two open intervals $\mathbb{X}^{\mathrm{C}}_i\subset \mathbb{X}^{\mathrm{B}}_i\subset{\mathbb{X}}_i$  and  define the following structure:
\begin{itemize}[nosep]
\item {\em base}: $\Pi^{\mathrm{B}}_i
=\Pi_i\cap(
{\mathbb{X}}^{\mathrm{B}}_i\times\mathbb Y_i\times\mathbb Z_i),
\quad H^{\mathrm{B}}=H\cap \Pi^{\mathrm{B}}_i,
\quad V^{\mathrm{B}}=g(\tilde H^{\mathrm{B}})\cap \Pi^{\mathrm{B}}_i,
$
\item {\em center}: $\Pi^{\mathrm{C}}_i=\Pi_i\cap(
{\mathbb{X}}^{\mathrm{C}}_i\times\mathbb Y_i\times\mathbb Z_i),
\quad H^{\mathrm{C}}=H\cap \Pi^{\mathrm{C}}_i,
\quad V^{\mathrm{C}}=g(\tilde H^{\mathrm{C}})\cap \Pi^{\mathrm{B}}_i$,
\end{itemize}
where $H$ and $V$ are any horizontal and vertical strips of $ \Pi_i$, with $V=g(\tilde  H)\cap \Pi_i$ for some $\tilde H$ of $\Pi_j$.
We say that an ss-disc $S$ {\em crosses} $\Pi_i^{\B/\C}$ if $S\cap \Pi_i\subset \Pi_i^{\B/\C}$; similarly for the strips.
\end{defi}

We say that two sets are {\em separated} if there is no ss-disc intersecting both of them. 
\begin{defi}[Base-Center  covering property]\label{defi:BC}
 A cs-Markov partition $(\{\Pi_i\},g)$ has the {\em Base-Center covering property}  if, in addition to~\ref{word:A1}, the cubes $\Pi_i$ can be structured  such that, with $S\subset  \Pi^\B_i$ and $V\subset \Pi_i$ denoting ss-discs and vertical strips, the following are satisfied  for every $i$:
\begin{itemize}[nosep]
\item[\setword{(A2)}{word:A2}] if $S\cap \partial_x V_1^{\mathrm{B}}\neq \emptyset$, then  $S$ crosses some center $V^\mathrm{C}_2$ which
is {separated} from $V^{\mathrm{C}}_1$, and
\item[\setword{(A3)}{word:A3}] if $S\cap {\Pi}_i^{\mathrm{C}}\neq \emptyset$, then $S$  crosses some center  $V^{\mathrm{C}}$.
\end{itemize}
\end{defi}

\begin{figure}[H]
\begin{center}
\includegraphics[scale=1.2]{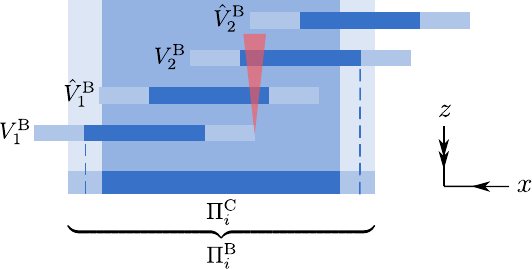}
\caption{The fulfilment of properties~\ref{word:A1}--\ref{word:A3}, where the red  area is covered by the strong-stable cone field $\mathcal{C}^{ss}$ }
\label{fig:BC}
\end{center}
\end{figure}

Property~\ref{word:A2} says that the vertical strips are `aligned', while property~\ref{word:A3} requires sufficiently many vertical strips whose centers cover  ${\Pi}_i^{\mathrm{C}}$. Figure~\ref{fig:BC} provides an illustration of properties~\ref{word:A1}--\ref{word:A3}, where \ref{word:A2} for the strips $V_1$ and $V_2$ are drawn explicitly with cone field $\mathcal{C}^{ss}$ (in red).

\begin{rem}\label{rem:coreblender}
Denote 
$$\mathcal{W}=\{\ell^u \subset W^u(\Lambda): g^{-n}(\ell^u) 
\mbox{ belongs to some }\Pi^{\mathrm{C}}_{i_n} \mbox{ for every } n\geq 0
\}.$$
Condition \ref{word:A3} implies that any ss-disc crossing $\Pi^{\mathrm{C}}_i$ intersects some unstable leaf $\ell^u\in \mathcal{W}$. To see this, one just needs to repeat the proof of Proposition~\ref{prop:nontr}, with $\Pi^{\mathrm{C}}_{i_n}$ in place of $\Pi_{i_n}$.
\end{rem}


\begin{defi}[Separated standard blenders]\label{defi:ablender}
A standard blender is called {\em separated} if its inducing Markov partition  satisfies the aligned covering property.
\end{defi}

One can check that blender-horseshoes are separated  cu-blenders induced from a cu-Markov partition with two elements. 

\subsection{Creation of a robust homoclinic tangency. Proof of Theorem~\ref{thm:tangency}}\label{sec:equi}


Robust tangencies  were achieved in \citep[Theorem 4.9]{BonDia:12b} for blender-horseshoes. Here we prove it for all separated  blenders. The proof is based on the notion of {\em folding manifolds} introduced in \citep{BonDia:12b}, which is defined below with modifications to fit the general setting of the present paper. We will see that Theorem~\ref{thm:tangency} is a direct consequence of the fact {\em every folding manifold has a tangency with $W^u(\Lambda)$}. The proof of this fact contains two steps. First, we  show in Lemma~\ref{lem:folding2} that if a folding manifold has no tangency with $W^u(\Lambda)$, then  its preimage must contain another folding manifold. Next, we show  in Proposition~\ref{prop:folding} that such property  implies the existence of  a tangency.


Let  $\Lambda$ be a separated  cs-blender with the inducing Markov partition $(\{\Pi_i\},g)$, and let $\mathcal{W}$ be the collection of unstable leaves defined in Remark~\ref{rem:coreblender}.

\begin{defi}[Folding manifolds]\label{defi:folding}
A submanifold $\mathcal{S}\subset \Pi^\B_i$ of dimension $(d_{ss}+1)$ is called a {\em folding manifold} of the blender $\Lambda$ if there exists a smooth family $\{S_t\}_{t\in[0,1]}$ of  ss-discs crossing $\Pi^\B_i$ such that
\begin{itemize}[nosep]
\item $\mathcal{S}=\bigcup_{t\in[0,1]} S_t$, and
\item $S_0$ and $S_1$ intersect a leaf $\ell^u \in \mathcal{W}$.
\end{itemize}
\end{defi}


\begin{figure}[h]
\begin{center}
\includegraphics[scale=.6]{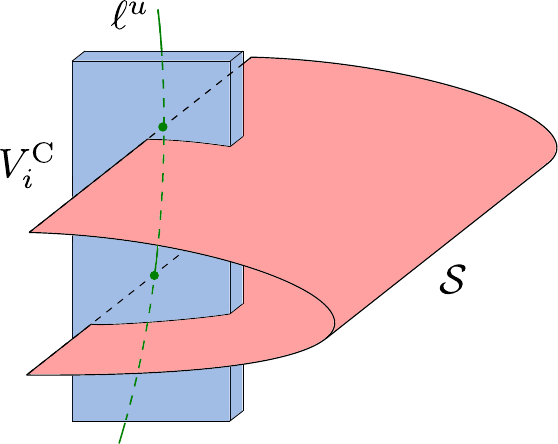}
\caption{A folding manifold}
\label{fig:folding1}
\end{center}
\end{figure}

See Figure~\ref{fig:folding1} for an illustration. 

\begin{rem}\label{rem:differ_folding}
The definition of folding manifolds in \citep{BonDia:12b} requires $\mathcal{S}$  to intersect the local unstable manifolds of two prescribed saddles in the two vertical strips that constitute the Markov partition of a blender-horseshoe. In contrast, here we only need $\mathcal{S}$ to attach to an arbitrary leaf in $\mathcal{W}$.
\end{rem}

For any u-disc $S^u$ with a defining function $(x,z)=(s^u_1(y),s^u_2(y))$, consider the family $\{S^u_t\}$ of u-discs with each $S^u_t$  given by $(x,z)=(s^u_1(y)+t,s^u_2(y))$. Take any  interval $I\subset\mathbb{R}$  such that $\bigcup_{t\in {I}}S^u_t$ covers $\Pi_i$. We denote 
\begin{equation}\label{eq:fakeleaf}
\mathcal{L}(S^u)=\Pi_i\cap\bigcup_{t\in{I}}S^u_t.
\end{equation}


\begin{lem}\label{lem:folding1}
Any folding manifold $\mathcal{S}$ has some point $P\in \mathcal{S}$ such that $T_P\mathcal{S}\cap \mathcal{C}^u\neq \emptyset$.
\end{lem}
\begin{proof}
 Let $\ell^u$ be the leaf  defining the folding manifold, and $\mathcal{L}:=\mathcal{L}(\ell^u)$ be defined as in \eqref{eq:fakeleaf}. Denote by $P_0=(x_0,y_0,z_0)$ and $P_1=(x_1,y_1,z_1)$ the intersections points of $S_0\cap \ell^u$ and $S_1\cap \ell^u$, respectively. Take the smooth curve $\gamma\subset \mathcal{S}\cap\mathcal{L}$ connecting $P_0$ and $P_1$ such that $\gamma(t)=(\gamma_x(t),\gamma_y(t),\gamma_z(t))\in S_t$ with $\gamma(0)=P_0$ and $\gamma(1)=P_1$. 
 
Let the defining function of $\ell^u$ be  $(x,z)=(\s^u_1(y),\s^u_2(y))$. By construction we have 
$\gamma_z(t)=\s^u_2(\gamma_y(t))$. Since $\ell^u$ is tangent to $\mathcal{C}^u$, we have
$$
\|\gamma'_z(t)\|=\left\|\dfrac{d \s^u_2(\gamma_y(t))}{dy}\cdot \gamma'_y(t)\right\|<K^u\cdot \|\gamma'_y(t)\|,
 $$
where $K^u$ is the cone constant in \eqref{eq:cone}. It follows that,  if $\gamma'_x(t)$ vanishes at some $t'\in(0,1)$, then   $T_{\gamma(t')}\mathcal{S}\cap \mathcal{C}^u\neq \emptyset$. Now suppose that 
$T_{\gamma(t)}\mathcal{S}\cap \mathcal{C}^u= \emptyset$ for all $t$. In particular, $\gamma'_x(t)$ never vanishes, and hence  $|\gamma_x(t)|$ is strictly increasing. Since the tangent spaces of $\gamma$ are transverse to $\mathcal{C}^u$ everywhere, for any given variation in $y$-coordinates, say $\|y_1-y_2\|$, the corresponding variation $\Delta x_\gamma$ in the $x$-coordinate of $\gamma$ must be strictly larger than that of any curve tangent to $\mathcal{C}^u$. However, by construction, we have $\Delta x_\gamma=|x_1-x_0|$, which is equal to the $x$-variation of any smooth curve in $\ell^u$ connecting $P_0$ and $P_1$.
\end{proof}

\begin{lem}\label{lem:folding2}
If $\mathcal{S}$ is a folding manifold of  $\Lambda$, then  either $\mathcal{S}$ has a non-transverse intersection with $W^u(\Lambda)$, or the preimage $g^{-1}(\mathcal{S})$  contains a folding manifold.
\end{lem}

\begin{proof}
Let $V_0$ be the vertical strip containing the leaf $\ell^u$ of the folding manifold.

\noindent (Case 1): If every ss-disc in $\mathcal{S}$ crosses  $V_0$, then $g^{-1}(\mathcal S\cap V_0)$ is a folding manifold. Indeed, by ~\ref{word:A1} this preimage lies in some $\Pi^\B_j$ and consists of ss-discs, and  the leaf $g^{-1}(\ell^u)$ belongs to $\mathcal{W}$ by definition. The projection of this case to the $(x,y)$-space is sketched in the left picture of Figure~\ref{fig:folding2}.

\begin{figure}[!h]
\begin{center}
\includegraphics[scale=.7]{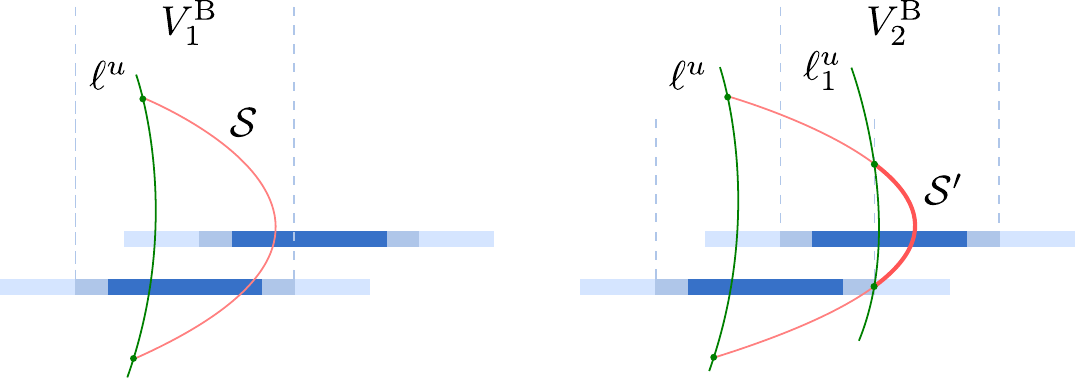}
\caption{Two cases of the position of a folding manifold}
\label{fig:folding2}
\end{center}
\end{figure}

\noindent (Case 2): Now suppose that  $V_0$ is not crossed by all ss-discs of $\mathcal{S}$. In particular, there exists  $S_{t_1}$ intersecting $\partial_x V_0^\mathrm{B}$ for some $t_1\in (0,1)$. By ~\ref{word:A2} there is some vertical strip $V_1$ such that $S_{t_1}$ crosses $V_1^\mathrm{C}$. 
It  follows from Remark~\ref{rem:coreblender} that  $S_{t_1}$ intersects a leaf $\ell^u_1\in \mathcal{W}$. If $\ell^u_1$ is tangent to $\mathcal{S}$, then we are done.
\begin{claim*}
If $\ell^u_1$ is  transverse to $\mathcal{S}$, then $\ell^u_1$ intersects $S_{t_2}$ for some $t_2\in(0,1)\setminus\{t_1\}$.
\end{claim*}
Let us prove the claim. Let $\ell^u$ be the leaf in the definition of the folding manifold. We attach to each point of $\ell^u$  an ss-disc such that they form a continuous family  including $S_0$ and $S_1$. Call its union $\mathcal{H}$. Next, consider the cu-surface $\mathcal{L}:=\mathcal{L}(\ell^u_1)$ as in \eqref{eq:fakeleaf}. By dimensional reason, $\ell^u_1$ divides $\mathcal{L}$ into  two connected components $C_\L$ and $C_\R$, the left and right ones in the obvious sense. The  curve $\mathcal{H}\cap\mathcal{L} $ must lie entirely in either $C_\L$ or $C_\R$, for otherwise there exists an ss-disc intersecting both $\ell^u\subset V_0^\mathrm{C}$ and $\ell^u_1\subset V_1^\mathrm{C}$, contradicting~\ref{word:A2}. In particular, we have that $S_0\cap \mathcal{L}$ and $S_1\cap \mathcal{L}$ lie in the same component. The claim then follows from the transversality assumption as it says that, for a small neighborhood $I$ of $t_1$, the intersection $(\bigcup_{t\in I}S_t)\cap\mathcal{L}$ has points in both components.

By construction, $\mathcal{S}':=\bigcup_{t\in[t_1,t_2]}S_t$ is a folding manifold. See the right picture of Figure~\ref{fig:folding2} for an illustration.  If every ss-disc of $\mathcal{S}'$ crosses the base of $V_1$, then we go to case 1 and finish the proof. Otherwise, we repeat the above procedure and obviously we must arrive at  case 1 after finitely many repetitions.
\end{proof}


\begin{prop}\label{prop:folding}
If $\mathcal{S}$ is a folding manifold of a separated standard cs-blender $\Lambda$, then $\mathcal{S}$ has a non-transverse intersection with $W^u(\Lambda)$.
\end{prop}
The proof of this proposition is the same as that of \citep[Proposition 4.4]{BonDia:12b}, where one replaces \citep[Remark 4.3]{BonDia:12b} and \citep[Lemma 4.5]{BonDia:12b} by Lemmas~\ref{lem:folding1} and~\ref{lem:folding2}, respectively.  For completeness, we reprove it here.

\begin{proof}[Sketch of proof.]
The  idea is similar to that of  the proof of Proposition~\ref{prop:nontr}, with a folding manifold $\mathcal{S}$ in place of an ss-disc. Suppose in Lemma~\ref{lem:folding2} there are no tangencies found in every backward iterate of the folding manifold $\mathcal{S}$. Then, we obtain a nested sequence $\{\hat {\mathcal{S}}_n\subset \mathcal{S}\}_{n\geq 0}$ with $\mathcal{S}_0=\mathcal{S}$ such that $\mathcal{S}_n:=g^{-n}(\hat{\mathcal{S}}_n)$ is a folding manifold for each $n$, and $\mathcal{S}':=\bigcap \hat {\mathcal{S}}_n\subset W^u(Q)$ for some point $Q\in\Lambda$ by the connectedness of $\mathcal{S}'$.

Lemma~\ref{lem:folding1} implies that there exists for each $n$ a point $P_n\in \mathcal{S}_n$ such that $T_{P_n} \mathcal{S}_n$ contains a non-zero vector $v_n$ which also lies in $ \mathcal{C}^u$. Denote $\hat P_n=g^n(P)\in \hat{\mathcal{S}}_n$ and $\hat v_n=Dg^n(v_n)\in T_{\hat P_n}\hat{\mathcal{S}}_n$.   Passing to subsequences if necessary, one has 
$$
\hat P_n\to \hat P\in \mathcal{S}' 
\quad\mbox{and}\quad
\hat v_n\to \hat v \in T_{\hat P}\mathcal{S}'
$$
for some point $\hat P$ and vector $\hat v$. Let us denote by $\mathcal{C}^u_P$ the cone of cone field $\mathcal{C}^u$ at point $P$. 
By the uniform hyperbolicity,  the size of the iterated cone $Dg^n(\mathcal{C}^u_P)$ goes to zero exponentially as $n\to +\infty$. Since one has the hyperbolic splitting $T_P\Pi_i=E_P^{ss}\oplus E_P^{cs}\oplus E^u_P$, it follows that $Dg^n(\mathcal{C}^u_{P_n})$, and hence $\hat v_n$, converge to the unstable space $E^u_{\hat P_n}$ as $n\to +\infty$. Thus, by continuity of the splitting one has $\hat v\in E^u_{\hat P}$, which is just $T_{\hat P}W^u(Q)$.
\end{proof}

\begin{proof}[Proof of Theorem~\ref{thm:tangency}]
We prove   only for a separated  cs-blender $\Lambda$, and the proof for a cu-blender is the same. By construction, for any local unstable leaf $\ell^u\in\Pi_i$, the union of its iterates  $\bigcup g^n(\ell^u)$ is dense in $W^u(\Lambda)$. Let us further take $\ell^u$ from $\mathcal{W}$ (defined in Remark \ref{rem:coreblender}).
By an arbitrarily $C^r$-small perturbation, we can find a piece $W^s$ of $W^s(\Lambda)$ tangent to $\ell^u$. With an additional small perturbation if necessary, we make the tangency quadratic. Then, since $W^s$ is foliated by strong-stable leaves which  particularly are ss-disc, unfolding the tangency will create inside $W^s$ a folding manifold attached to $\ell^u$.  The theorem now follows from Proposition~\ref{prop:folding} and the fact that  being a folding manifold is a $C^1$-robust property.
\end{proof}

\section{Arrayed standard blenders}\label{sec:array}
We  now introduce  arrayed  blenders. Roughly speaking, a $k$-arrayed  blender is a proper combination of $k$  standard blenders such that the central direction is covered $k$ times in an aligned fashion.

\subsection{Base-Center-Gap structure}
 
\begin{figure}[h]
\begin{center}
\includegraphics[scale=.6]{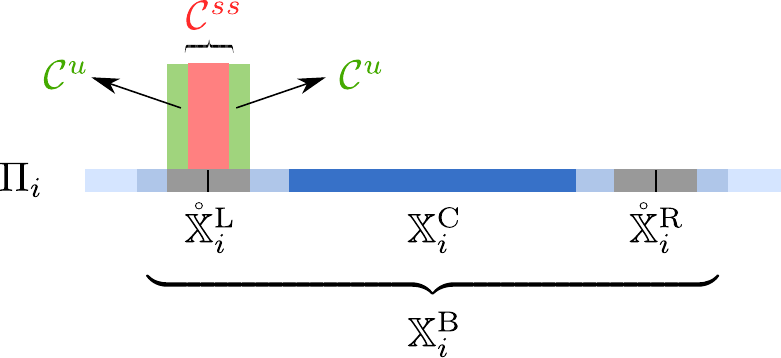}
\caption{The Base-Center-Gap structure projected to $\X_i$}
\label{fig:BCG}
\end{center}
\end{figure}

Recall that we say two sets $A$ and $B$ are separated if  no ss-disc intersecting both of them. We now denote this by $A \parallel B$, and write $A \nparallel B$ if there is some ss-disc connecting them.
By the partial hyperbolicity condition in Definition~\ref{defi:markov}, the first component of the map $g|_H$ for a horizontal strip $H$ is strictly monotone in $x$. We say that $g|_H$ is orientation-preserving if this component is increasing in $x$, and orientation-reversing if it is decreasing in $x$.
\begin{defi}[Base-Center-Gap (BCG) structure of $\Pi_i$]\label{defi:structure2}
For each $i$ consider two open intervals $\mathbb{X}^{\mathrm{C}}_i$ and $\mathbb{X}^{\mathrm{B}}_i$ satisfying 
  $\overline{\mathbb{X}}^{\mathrm{C}}_i\subset\mathbb{X}^{\mathrm{B}}_i\subset\mathbb{X}_i$, and two points $x^\mathrm{L}_i,x^\mathrm{R}_i$ in the left and, respectively, right (by the usual orientation on $\mathbb{R}$) components of $\mathbb{X}^{\mathrm{B}}_i\setminus \overline{\mathbb{X}}^{\mathrm{C}}_i$, and denote
  $$
  \mathring{\Pi}^{\mathrm{L}}_i:=\Pi_i\cap(\{x^\mathrm{L}_i\}\times\mathbb Y_i\times\mathbb Z_i)
  \quad\mbox{and}\quad
  \mathring{\Pi}^{\mathrm{R}}_i:=\Pi_i\cap(\{x^\mathrm{R}_i\}\times\mathbb Y_i\times\mathbb Z_i).  
  $$
We define
\begin{itemize}[nosep]
\item {\em base}: $\Pi^{\mathrm{B}}_i
=\Pi_i\cap(
{\mathbb{X}}^{\mathrm{B}}_i\times\mathbb Y_i\times\mathbb Z_i),
$
\item {\em center}: $\Pi^{\mathrm{C}}_i
=\Pi_i\cap({\mathbb{X}}^{\mathrm{C}}_i\times\mathbb Y_i\times\mathbb Z_i)$,
\item {\em left gap}: $\Pi^\L_i$, the union of all u-discs $S^u$ satisfying $S^u\nparallel \mathring{\Pi}^{\mathrm{L}}_i,$
\item {\em right gap}: $\Pi^\R_i$, the union of all u-discs $S^u$ satisfying $S^u\nparallel \mathring{\Pi}^{\mathrm{R}}_i. $
\end{itemize}
For any horizontal strip $H$ and vertical strip $V$   with $V=g(\tilde  H)\cap \Pi_i$ for some $\tilde H$ of $\Pi_j$, we denote
$$
H^{\B/\C/\L/\R}=H\cap \Pi^{\B/\C/\L/\R}_i,\quad
\mathring{H}^{\L/\R}=H\cap \mathring{\Pi}^{\L/\R}_i,\quad
V^{\B/\C}=g(\tilde H^{\B/\C})\cap \Pi^{\B}_i,
$$
The gaps for vertical strips are defined as follows: if $g|_H$ is orientation-preserving, then
$$\mathring{V}^{\L/\R}=g(\mathring{H}^{\L/\R})\cap \Pi_i^\B,\quad
V^{\L/\R}=g(H^{\L/\R})\cap \Pi^\B_i;$$
if $g|_H$ is orientation-reversing, then
$$\mathring{V}^{\L/\R}=g(\mathring{H}^{\R/\L})\cap \Pi_i^\B,\quad
V^{\L/\R}=g(H^{\R/\L})\cap \Pi^\B_i;$$
We say that an ss-disc $S$ {\em crosses} $\Pi_i^{\B/\C/\L/\R}$  if $S\cap \Pi_i\subset \Pi_i^{\B/\C/\L/\R}$; similarly for the strips.
\end{defi}
See Figure~\ref{fig:BCG} for an illustration, where the size of the gaps (in gray) are determined by the size of the cone fields.

For a  finite collection $\mathcal{V}$ of different vertical strips in $\Pi_i$, its base, center and gaps, denote by $\mathcal{V}^{\mathrm{B}},\mathcal{V}^{\mathrm{C}},\mathcal{V}^{\mathrm{L}},\mathcal{V}^{\mathrm{R}}$, are defined as the collections of the corresponding parts of the vertical strips in it. 
\begin{defi}[Arrays of vertical strips]
A finite collection $\mathcal{V}$ of $k$ different vertical strips in $\Pi_i$  is called a {\em vertical $k$-array} of $\Pi_i$ if there exists some ss-disc crossing every $V\in \mathcal{V}$; such ss-discs are said to {\em cross the array}, and the same for crossing of base, center, gaps of $\mathcal{V}$. 
\end{defi}

Let $\partial_x^\mathrm{L/R} \Pi_i$ be the left/right $x$-boundary in the obvious sense (see the discussion above Definition~\ref{defi:disc}), and the same for $\partial_x^\mathrm{L/R} H$. Let $V\subset \Pi_i$ be of the form $V=g(H)\cap\Pi_i$ for some $H\subset\Pi_j$.  If $g|_H$ is orientation-preserving, then
$$\partial_x^\mathrm{L} V=g(\partial_x^\mathrm{L} H)\cap \Pi_i,\quad
\partial_x^\mathrm{R} V=g(\partial_x^\mathrm{R} H)\cap \Pi_i,$$
and, if $g|_H$ is orientation-reversing in $x$, then
$$\partial_x^\mathrm{L} V=g(\partial_x^\mathrm{R} H)\cap \Pi_i,\quad
\partial_x^\mathrm{R} V=g(\partial_x^\mathrm{L} H)\cap \Pi_i.$$

\begin{defi}[Base-Center-Gap (BCG)  covering property for arrays]\label{defi:bcgcover}
A cs-Markov partition $(\{\Pi_i\},g)$ has the {\em $k$-arrayed BCG covering property} if $\Pi_i$ can be structured  such that, with $S\subset \Pi^\B_i,V\subset \Pi_i,\mathcal{V}\subset \Pi_i$ denoting its ss-discs, vertical strips and, respectively,  vertical $k$-arrays, the following are satisfied:
\begin{itemize}[nosep]
\item[\setword{(B1)}{word:B1}] any $S$ crossing $\Pi^{\mathrm{B}}_i$  crosses at least one vertical  array $\mathcal{V}^\mathrm{B}$ of bases,

\item[\setword{(B2)}{word:B2}] if $S\cap \partial^{\mathrm{R}}_x V^{\mathrm{B}}\neq \emptyset$ for some $V\in \mathcal{V}_1$, then  $S$ crosses some $\mathcal{V}^\mathrm{C}_2$ such that $ \mathcal{V}_2^\mathrm{B} \parallel \mathcal{V}_1^\mathrm{L}$,
\item[\setword{(B2')}{word:B2'}] if $S\cap \partial^{\mathrm{L}}_x V^{\mathrm{B}}\neq \emptyset$ for some $V\in \mathcal{V}_1$, then  $S$ crosses some $\mathcal{V}^\mathrm{C}_2$ such that $ \mathcal{V}_2^\mathrm{B} \parallel \mathcal{V}_1^\mathrm{R}$,
\item[\setword{(B3)}{word:B3}] if $S\cap\partial^{\mathrm{L}}_x\Pi_i^{\mathrm{C}}\neq \emptyset$, then   $S$ crosses some  $\mathcal{V}^{\mathrm{C}}$ such that $ \mathcal{V}^\mathrm{B} \parallel \Pi_i^\mathrm{L}$,
\item[\setword{(B3')}{word:B3'}] if $S\cap\partial^{\mathrm{R}}_x\Pi_i^{\mathrm{C}}\neq \emptyset$, then   $S$ crosses some  $\mathcal{V}^{\mathrm{C}}$ such that $ \mathcal{V}^\mathrm{B} \parallel \Pi_i^\mathrm{R}$, and  
\item[\setword{(B4)}{word:B4}] $\Pi^{\L/\R}_i\parallel \Pi^\C_i$ for any $\Pi_i$ and $V^{\L/\R}\parallel V^\C,V^{\L/\R}\parallel \partial^{\L/\R}_x V^\B$ for  any  $V$.
\end{itemize}
We say that a cu-Markov partition satisfies the BCG covering property if its inverse satisfies it.
\end{defi}

\begin{figure}[h]
\begin{center}
\includegraphics[scale=.9]{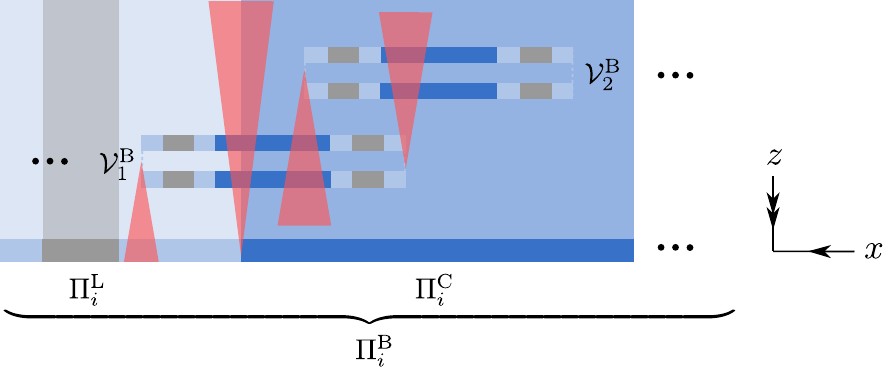}
\caption{The fulfilment of properties~\ref{word:B2}--\ref{word:B4}, where the red  area is covered by the strong-stable cone field $\mathcal{C}^{ss}$}
\label{fig:BCGcover}
\end{center}
\end{figure}

Property~\ref{word:B1} is just the covering property~\ref{word:A1}, but stated for arrays. Figure~\ref{fig:BCGcover} provides an illustration of the fulfilment of other properties, except for the symmetric ones:~\ref{word:B3'} and $\Pi^\R_i\parallel \Pi^\C_i$ in~\ref{word:B4}, where, in particular, properties~\ref{word:B2},~\ref{word:B2'} and~\ref{word:B3} are shown explicitly with the cone field $\mathcal{C}^{ss}$ (in red).



\begin{defi}[$k$-arrayed standard blenders]\label{lem:kblender}
Given a partially hyperbolic Markov partition that satisfies the $k$-arrayed BCG covering property, its   locally maximal set is called a {\em $k$-arrayed standard blender}. The blender is further called center-stable (cs) for cs-Markov partitions and center-unstable (cu) for cu-Markov partitions.
\end{defi}

\begin{rem}\label{rem:diff}
Comparing Definitions~\ref{defi:BC} and~\ref{defi:bcgcover}, one sees that an arrayed  blender is not necessarily a separated one, and vice versa. The reason is that the BCG covering property enables us to distinguish the left and right (to define the folding manifolds in Section~\ref{sec:reproducing}), based on the fact that the central dynamics in $\X$ is one dimensional. In contrast, the BC covering property essentially provides the notions of inside and outside of $\Pi^\C_i$ so that one can define folding manifolds  which are large enough (in the central direction) -- it crosses the boundary $\partial \X^\C_i$ and hence have parts both inside and outside $\Pi^\C_i$.  Since  $\partial\Pi_i^\C$  is essentially determined by $\partial \X^\C_i$,  the BC covering property can be generalized to the case of $d_c$-dimensional $(d_c>1)$ center by distinguish the inside and outside of $\partial \X^\C_i=S^{d_c}$, the $d_c$-dimensional sphere. Hence, we believe that  separated  blenders can be constructed for Markov partitions with higher center dimension to create robust tangencies of higher corank, that is, tangencies at which the two tangent spaces share a common subspace of  dimension larger than 1.  
\end{rem}

\subsection{Reproducing lemma}\label{sec:reproducing}

Let $\mathcal{S}\subset \Pi^\B_i$ be a submanifold of dimension $(d_{ss}+1)$  which is the union of a smooth family $\{S_t\}_{t\in[0,1]}$ of  ss-discs crossing $\Pi^\B_i$.

\begin{defi}[Prefolding  manifolds]\label{defi:folding2'}
The submanifold $\mathcal{S}$ is called a {\em right/left-prefolding manifold} if
\begin{itemize}[nosep]
\item $S_0$ and $S_1$ intersect a leaf $\ell^u\subset \left( W^u(\Lambda)\cap \Pi^{\mathrm{L/R}}_i\right)$, and 
\item $S_t\cap \Pi^\mathrm{C}_i\neq \emptyset$ for some $t\in (0,1)$.
\end{itemize}
\end{defi}

\begin{defi}[Folding manifolds]\label{defi:folding2}
The submanifold $\mathcal{S}$ is called a {\em right/left-folding manifold related to a vertical strip $V$} if \begin{itemize}[nosep]
\item $S_0$ and $S_1$ intersect a leaf $\ell^u\subset \left( W^u(\Lambda)\cap V^{\mathrm{L/R}}\right)$,  
\item $S_t\cap \ell^u=\emptyset$ for all $t\in (0,1)$, and 
\item $S_t\cap V^\mathrm{C}\neq \emptyset$ for some $t\in (0,1)$;
\end{itemize}
it is further called  {\em exact} if every ss-disc of $\mathcal{S}$ crosses $V^\B$.
\end{defi}

\begin{figure}[h]
\begin{center}
\includegraphics[scale=.9]{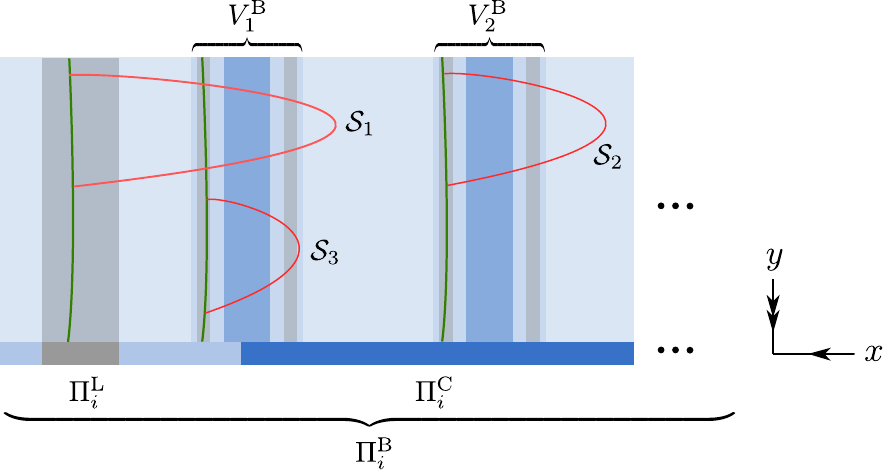}
\caption{A right-prefolding manifold ($\mathcal{S}_1$), a right-folding manifold ($\mathcal{S}_2)$ related to $V_2$, and an exact right-folding manifold  $(\mathcal{S}_3)$ related to $V_1$}
\label{fig:folding3}
\end{center}
\end{figure}

See Figure~\ref{fig:folding3} for an illustration.

\begin{defi}[Folding arrays]\label{defi:folding3}
The submanifold  $\mathcal{S}$ is called a {\em right/left-folding  array related to a  vertical $k$-array $\mathcal{V}=\{V_i\}$} if, for  $i=1,\dots,k$, there exist  unstable leaves $\ell^u_i$ and $t_i,t'_i\in(0,1]$ with $t_i<t'_i$ such that $\mathcal{S}_i:=\bigcup_{t\in[t_i,t'_i]} S_t$ are right-/left-folding manifolds related to $V_i$. It is further called  {\em exact} if all $\mathcal{S}_i$ are  exact.
\end{defi}

\begin{lem}[Reproducing lemma]\label{lem:re}
Let the $k$-arrayed BCG covering property be satisfied. If $\mathcal{S}$ is an exact folding array related to $\mathcal{V}=\{V_i\}^k_{i=1}$, then for each $i$  the preimage $g^{-1}(\mathcal{S}\cap V_i)$  contains an exact folding array in some $\Pi_{j}$. 
\end{lem}

The novelty in this result is the nonstop reproducing process, in contrast with Lemma~\ref{lem:folding2} where the process terminates upon detecting a  tangency (since  the arguments there for generating a new folding manifold is  based on the absence of tangencies). This property is crucial for obtaining uncountably many tangencies in Proposition~\ref{prop:foldarray}. 

Since for each $i$ the intersection $\mathcal{S}\cap V_i$ contains an exact folding manifold related to $V_i$, Lemma~\ref{lem:re} is an immediate consequence of the following three lemmas.

\begin{lem}\label{fact:1}
The preimage by $g$ of an exact folding manifold is a prefolding manifold.
\end{lem}

\begin{lem}\label{fact:2}
Every prefolding manifold contains a folding array.
\end{lem}

\begin{lem}\label{fact:3}
Every folding array contains an exact folding array.
\end{lem}

\begin{proof}[Proof of Lemma~\ref{fact:1}]
Let $\mathcal{S}$ be an exact folding manifold related to ${V}$. It follows from definitions of the BCG structure and folding manifolds that one has $g^{-1}(\mathcal{S}\cap V^{\mathrm{L/R/C}})\subset \Pi_j^{\mathrm{L/R/C}}$ for some $j$ if $g^{-1}|_V$ is orientation preserving or $g^{-1}(\mathcal{S}\cap V^{\mathrm{L/R/C}})\subset \Pi_j^{\mathrm{R/L/C}}$ if $g^{-1}|_V$ is orientation reversing. Since the cone field $\mathcal{C}^{ss}$ is backward-invariant, in either case  the preimage is a prefolding manifold.
\end{proof}

\begin{proof}[Proof of Lemma~\ref{fact:2}]
Let $\mathcal{S}=\bigcup_{t\in[0,1]}S_t$ be a right-prefolding manifold of $\Pi_i$. Since some $S_t$ crosses $\Pi_i^\C$ by definition, there exists $S_{t_0}$ such that $S_{t_0}\cap\partial^\mathrm{L}_x\Pi^\mathrm{C}_i\neq \emptyset$. Let $\mathcal{V}$ be the vertical $k$-array given by~\ref{word:B3}.  Take any $V_j\in\mathcal{V}$. Since $S_{t_0}$ crosses $\mathcal{V}^\mathrm{C}$, it in particular crosses $V^\mathrm{C}_j$. We prove the lemma in several steps.

\noindent\textit{Step 1: there exists $t_1\in(0,t_0)$ such that $S_{t_1}\cap \mathring{V}_j^\mathrm{L}\neq\emptyset$ (see Definition~\ref{defi:bcgcover}).}

By the continuity of $S_t$ on $t$ and the fact that ss-discs
are defined for all $z$ values (see Definition~\ref{defi:disc}),  there must exist some $t_1\in(0,t_0)$ such that  $S_{t_1}\cap \mathring{V}_j^\mathrm{L}\neq\emptyset$ or $S_{t_1}\cap \mathring{V}_j^\mathrm{R}\neq\emptyset$. We prove below that  the latter case implies the existence of some $t'\in(0,t_0)$ with $S_{t'}\cap \mathring{V}_j^\mathrm{L}\neq\emptyset$. Denote by $\ell^u$ the unstable leaf in the definition of a prefolding manifold. Let us attach to each point of $\ell^u$ an ss-disc, including $S_0$ and $S_1$, to form a su-surface $\mathcal{H}$ of codimension one which divides $\Pi_i$ into two connected components, the left one $C_\L$ and the right one $C_\R$ in the obvious sense. So, we have $\Pi_i\setminus \mathcal{H}=C_\L\bigcup C_\R$. Similarly, we take a family of ss-discs including $S_{t_0}$ such that every disc intersects $V_j^\C$. It forms another codimension-one surface $\mathcal{H}_0$, which lies in $C_\R$ by~\ref{word:B3} and further divides  $C_\R$ into two components $C_{\R_1}$ and $C_{\R_2}$, labeled from left to right. Property~\ref{word:B4} implies $V_j^\L\subset \Pi_i\setminus C_{\R_2}$ and $V_j^\R\subset C_{\R_2}$. Since~\ref{word:B3} also implies $V_j^\L\subset C_\R$, one has $V_j^\L\subset C_{\R_1}$.  The desired statement then follows immediately from the continuity  since before intersecting $\mathring{V}_j^\mathrm{R}$ the family $\{S_t\}$ must go through $C_{\R_1}$, and hence intersect $\mathring{V}_j^\mathrm{L}$.

\noindent\textit{Step 2:  $S_{t_1}$  intersects an unstable leaf $\ell^u_1\in W^u(\Lambda)\cap V_j^\L$.}

Property~\ref{word:B1} implies that every ss-disc crossing $\Pi^\B_i$ intersects $W^u(\Lambda)$ (see Proposition~\ref{prop:nontr}). In particular, there exists some unstable leaf $\ell^u_1\subset W^u(\Lambda)$ intersecting $S_{t_1}$, which lies in $V_j^\B$ by Remark~\ref{rem:unstableleaf}.
We in fact have $\ell^u_1\subset V_j^\L$. To see this, let $V_j=g(H_j)$ for some horizontal strip $H_j$. Assume that $g|_{H_j}$ is orientation preserving in $x$, so $\mathring{V}_j^\mathrm{L}=g(\mathring{H}_j^\mathrm{L})$.  By the invariance of the unstable foliation and the cone field $\mathcal{C}^{ss}$, we have that the ss-disc $(g|_{H_j})^{-1}(S_{t_1})$ intersects both $\mathring{H}_j^\mathrm{L}$ and the unstable leaf
$(g|_{H_j})^{-1}(\ell^u_1)$. Hence,
 $\ell^u_1\subset {V}_j^\mathrm{L}$  by the definition of gaps.  Obviously, the result also holds when $g|_{H}$ is orientation reversing in $x$. 
 
\noindent\textit{Step 3: there exists $t'_1\in [t_1,t_0)$ such that $S_{t'_1}$  intersects $\ell^u_1$ and $S_t\cap \ell^u_1= \emptyset$ for $t\in(t'_1,t_0]$.} 

Such $t'_1$ does not exist only if $S_t$ intersects $\ell^u_1$ for every $t\in [t_1,t_0]$. This is impossible since $S_{t_0}$ belongs to $\mathcal{H}_0$, the right boundary of $C_{\R_1}$, while $\ell_1^u\subset \overline V_j^\L\subsetneq C_{\R_1}$ by  \ref{word:B3} and \ref{word:B4}.

%

\noindent\textit{Step 4: there exists $t_2\in(t_0,1)$ such that $S_{t_2}$  intersects  $\ell^u_1$ and $S_t\cap \ell^u_1= \emptyset$ for $t\in[t_0,t_2)$.}

Since $S_{t_0}$ and  $S_{t_1}$ belong to the right, and, respectively, left boundaries of $C_{\R_1}$ while   $\ell_1^u\subset \overline V_j^\L\subsetneq C_{\R_1}$, the desired value $t_2$ follows from the continuity of $\{S_t\}$.

\noindent\textit{Step 5: completion of the proof.} 

 By construction, $\mathcal{S}_j=\bigcup_{t\in[t'_1,t_2]} S_t$ is a right-folding manifold of $V_j$. 
Repeat the above arguments for all $k$ vertical strips in $\mathcal{V}$, we find $k$ right-folding manifolds $\mathcal{S}_j$ and their union is a right-folding array related to $\mathcal{V}$. The proof for the case where $\mathcal{S}$ is a left-prefolding manifold is the same, with using~\ref{word:B3'} instead of~\ref{word:B3}. 
\end{proof}
 

\begin{rem}\label{rem:folding}
It can be seen from the proof that, with using \ref{word:B2} instead of \ref{word:B3}, if $\mathcal{S}$ is a right-folding manifold $\mathcal{S}$ related to some vertical strip $V$ (instead of a right-prefolding manifold) such that $\mathcal{S}\cap\partial^\R_x V\neq \emptyset$, then $\mathcal{S}$ contains a folding array given by \ref{word:B2}. Similarly, for left-folding manifolds,  using \ref{word:B2'}.
\end{rem}

\begin{proof}[Proof of Lemma~\ref{fact:3}]
Let $\mathcal{S}=\bigcup_{t\in[0,1]}S_t$ be a right-folding array of a vertical $k$-array $\mathcal{V}_0$. If $\mathcal{S}$ is not exact, then by definition there exists some $V_0\in\mathcal{V}_0$ such that the right-folding manifold $\mathcal{S}_0$ associated with $V_0$ is not exact. So, $\mathcal{S}_0$ contains some ss-disc $S_{t_0}$ such that $S_{t_0}$ intersects  $\partial^\L_x V_0^\mathrm{B}$ or $\partial^\R_x V_0^\mathrm{B}$.  

We claim that
$S_{t_0}\cap \partial^\R_x V_0^\B\neq \emptyset$. To prove the claim, first note that the gap $V_0^\L$ divides $V_0$ into two connected components $C_\L$ and $C_\R$, the left and right ones, such that $\partial^\L_x V_0^\B\subset C_\L$.  It is evident that if there exists an ss-disc intersecting both $\partial^\L_x V_0^\B$ and $ C_\R$, then one can construct another ss-disc intersecting $\partial^\L_x V_0^\B$ and $V^\L$. On the other hand, if  $\ell^u\subset V^\L_0$ is the unstable leaf associated with $\mathcal{S}_0$ as in Definition~\ref{defi:folding2}, then $\ell^u\cap S_{t_0}=\emptyset$, and hence $S_{t_0}$ intersects $ (V^\L \cup C_\R)$. The claim then follows, since $S_{t_0}\cap \partial^\L_x V_0^\B\neq \emptyset$ would contradict $\partial^\L_x V_0^\B\parallel V^\L$ given by ~\ref{word:B4}.

 Thus, by~\ref{word:B2}, there exists a vertical $k$-array $\mathcal{V}_1$ such that $S_{t_0}$ crosses $\mathcal{V}^\mathrm{C}_1$. 
 One then obtains a right-folding array of $\mathcal{V}_1$ by the same arguments in
the proof of Lemma~\ref{fact:2}, with $\mathcal{S}_0,V_0,\mathcal{V}_1$,~\ref{word:B2} in place of $\mathcal{S},\Pi_i,\mathcal{V}$,~\ref{word:B3}, see Remark \ref{rem:folding}.
 If this folding array is not exact, then we repeat this procedure. Since there are only finitely many vertical arrays and the new folding array obtained at each step is strictly smaller than the previous one, we will obtain an exact folding array after finite steps. The proof for the case where $\mathcal{S}$ is a left-folding array is the same, with  ~\ref{word:B2'} in place of~\ref{word:B2}.
\end{proof}

\subsection{Robust presence of uncountably many homoclinic tangencies. Proof of Theorem~\ref{thm:uncount}}\label{sec:robuncon}
Let $\Lambda$ be a $k$-arrayed   cs-blender induced from a cs-Markov partition $(\{\Pi_i\},g)$. Let us denote by $\mathcal{V}_{ij}$ the vertical arrays of $\Pi_i$ involved in the $k$-arrayed BCG covering property in Definition~\ref{defi:bcgcover}, where $j=1,\dots,j^*$ for some integer $j^* > 1$ depending on $i$. We label the corresponding vertical strips by $V_{ijs}\in \mathcal{V}_{ij}$ with $s=1,\dots,k$. We call the connected components of $W^u(\Lambda)\cap \Pi_i$ the local unstable leaves of $W^u(\Lambda)$ (in $\Pi_i$). Denote by $\Sigma^k_+:=\{1,\dots,k\}^{\mathbb{N}_0}$ the set of one-sided sequences of $k$ symbols.

\begin{prop}\label{prop:foldarray}
For any exact folding array $\mathcal{S}$ of  $\Lambda$, each  $\underline{s}=(s_n)_{n\geq 0}\in\Sigma^k_+$ corresponds to a sequence $(i_n,j_n,s_n)_{n\geq 0}$ such that the local unstable leaf defined by 
\begin{equation}\label{eq:tangencyleaf}
\ell^u(\underline{s})=\left(\bigcap_{n\geq 0} g^{n}(V_{i_nj_ns_n})\right)\cap \Pi_{i_0}
\end{equation}
 intersects  $\mathcal{S}$ non-transversely. Moreover, this correspondence is injective.
\end{prop}

\begin{proof}
Suppose $\mathcal{S}$ is related to $\mathcal{V}_{i_0j_0}$ in $\Pi_{i_0}$. By the exactness of $\mathcal{S}$, for any $s_0\in \{1,\dots,k\}$, the intersection $\mathcal{S}\cap V_{i_0j_0s_0}$ contains an exact folding manifold related to  $V_{i_0j_0s_0}$. The reproducing lemma (Lemma~\ref{lem:re}) implies that  the preimage $g^{-1}(\mathcal{S}\cap V_{i_0j_0s_0})$ contains an exact folding array related to an array $\mathcal{V}_{i_1j_1}$ in some $\Pi_{i_1}$, where $i_1$ and $j_1$ depend on the choice of $s_0$. 
For any $s_1\in \{1,\dots,k\}$, the intersection of the new exact folding array with $V_{i_1j_1s_1}$ contains an exact folding manifold. Invoking the reproducing lemma on $V_{i_1j_1s_1}$, we find another exact folding array related to some array   $\mathcal{V}_{i_2j_2}$. Repeating this  procedure infinitely many times, we obtain a sequence $\underline s=(s_n)_{n\geq 0} \in \Sigma^k_+$. By construction, it corresponds to a sequence $(i_n,j_n,s_n)_{n\geq 0}$ with $i_n$ and $j_n$ depending on $s_{n-1}$ for $n>0$  such that, with $\mathcal{S}_0(\underline{s}):=\mathcal{S}$,  each of the preimages $\mathcal{S}_{n+1}(\underline{s}):=g^{-1}(\mathcal{S}_{n}(\underline{s})\cap V_{i_n j_n s_n})$  contains an exact folding manifold.  Since Lemma~\ref{lem:folding1}  obviously holds for exact folding manifolds, one can argue as in the proof of Proposition~\ref{prop:folding} to conclude that $\mathcal{S}\cap W^u(\Lambda)$ has a non-transverse intersection point in $\bigcap_{n\geq 0} \hat{\mathcal{S}}_n(\underline{s})$, where $\hat{\mathcal{S}}_n(\underline{s}):=g^n(\mathcal{S}_n(\underline{s}))$. By construction, this intersection point belongs to the leaf defined by~\eqref{eq:tangencyleaf}, and  for different sequences the corresponding leaves are different.
\end{proof}

\begin{figure}[h]
\begin{center}
\includegraphics[scale=.9]{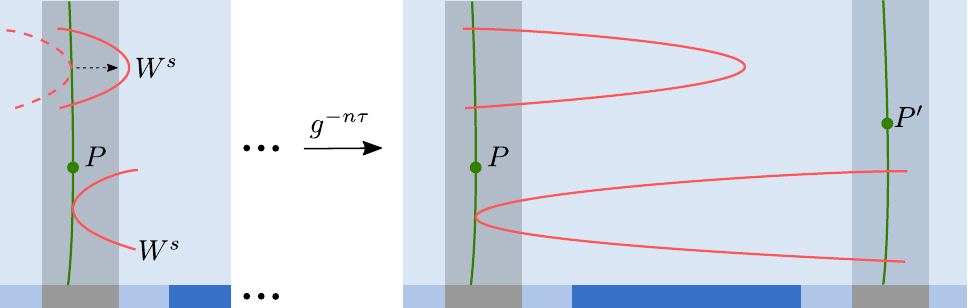}
\caption{Creating  prefolding manifolds in the two cases in the proof of Theorem \ref{thm:uncount}, projected to the $(x,y)$-space}
\label{fig:folding4}
\end{center}
\end{figure}

\begin{proof}[Proof of Theorem \ref{thm:uncount}]
By the symmetry of the problem, we only prove the result for a $k$-arrayed cs-blender $\Lambda$. The goal is to create an exact folding array in $W^s(\Lambda)$ and apply Proposition~\ref{prop:foldarray}. By Lemmas~\ref{fact:2} and~\ref{fact:3}, it  suffices to find a prefolding manifold. The robustness in the theorem is automatic since being a prefolding manifold or an (exact) folding array is a $C^1$-robust property.

Consider any $\Pi_i$. By \ref{word:B1} (and Proposition \ref{prop:nontr}), there exists a local unstable leaf $\ell^u\subset W^u(\Lambda)\cap \Pi^\L_i$. By construction  this leaf must contain a point in $\Lambda$. Indeed, $\ell^u$ intersects all local stable leaves in $\Pi_i$ and all intersection points belong to the locally maximal set $\Lambda$. Since periodic points are dense in $\Lambda$, we can take some periodic point $P$ such that $W^u_{loc}(P)\subset \Pi^\L_i$. Similarly, we find a periodic point $P'$ with $W^u_{loc}(P')\subset \Pi^\R_i$. 
Since $W^u(P)$ is dense in $W^u(\Lambda)$, by an arbitrarily small perturbation, we can create a quadratic tangency between a piece $W^s$ of $W^s(\Lambda)$ and $W^u_{loc}(P)$. There are now two cases.

The first case is when $W^s$ at the tangency point is pointing to the right in the $x$-direction, see the upper picture in Figure \ref{fig:folding4}.  In this case, we unfold the tangency by moving $W^s$ to the right in the $x$-direction and denote by $\mu$ this displacement. So the tangency exists at $\mu=0$, and there are two transverse intersection points of $W^s$ with $W^u_{loc}(P)$ at $\mu>0$. Note that $W^s$ is foliated by strong-stable leaves (which in particular are ss-discs), so these intersection points belong to two strong-stable leaves in $W^s$, denoted by $\ell^{ss}_0$ and $\ell^{ss}_1$.  Let $\tau$ be the period of $P$. Due to the contraction in $x$, for any fixed integer $n>0$, there exists $\mu_1>0$ such that $g^{-n\tau}(W^s)$ intersects the right boundary $\partial^\R_x \Pi_i$ at $\mu=\mu_1$. By continuity, one can find $\mu_2\in(0,\mu_1)$ such that $g^{-n\tau}(W^s)$ intersects $\Pi^\C_i$ but still lies in $\Pi^\B_i$ at $\mu=\mu_2$. Now, the part of $g^{-n\tau}(W^s)$  bounded by $g^{-n\tau}(\ell^{ss}_0)$ and $g^{-n\tau}(\ell^{ss}_1)$ is a right-prefolding manifold, as desired. Moreover, $\mu_2$, and hence the size of the perturbation, tend to zero as $n\to \infty$.

In the second case, $W^s$ at the tangency point is pointing to the left, see the lower picture in Figure \ref{fig:folding4}. We now only perturb the tangency to a quadratic one and do not unfold it. Let us write it as a union of strong-stable leaves, $W^s=\bigcup_{t\in[0,1]}\ell^{ss}_t$, where the tangency point belongs to $\ell^{ss}_{t_0}$ for some $t_0\in (0,1)$. Due to the contraction in $x$, the preimage $g^{-n\tau}(W^s)$ for a large $n$ must intersect $W^u_{loc}(P')$ in at least two leaves $g^{-n\tau}(\ell^{ss}_{t_1})$ and $g^{-n\tau}(\ell^{ss}_{t_2})$ satisfying $t_1\in (0,t_0)$ and  $t_1\in (t_0,1)$. Evidently,  the part of $g^{-n\tau}(W^s)$ bounded by these two leaves, i.e., $\bigcup_{t\in[t_1,t_2]}g^{-n\tau}(\ell^{ss}_{t})$, is a left-prefolding manifold.
\end{proof}

\section{Nearly-affine blender system (NABS)}\label{sec:hdc}
In this section, we prove Theorem~\ref{thm:hdc} by showing that the abstract constructions of separated and arrayed blenders in the previous sections arise naturally near heterodimensional cycles. Our starting point is the notion of a nearly-affine blender system, introduced below, which generalizes the renormalization results of   \citep{LiTur:24}.


Denote $d=\dim \mathcal{M}$.
Let $\{\Pi_\delta \}_{\delta>0}$ be a family of  embedded $d$-dimensional unit cubes, that is, there is a family of  charts $\{(U_\delta,\psi_\delta)\}_{\delta>0}$ with $\Pi_\delta\subset U_\delta$ and $\varphi_\delta(\Pi_\delta)=[-1,1]^d\subset \mathbb{R}^d$.  For each fixed $\delta$, we associate with $\Pi_\delta$ a family $\mathcal{F_\delta}:=\{F_{i,\delta}\}_{i\in\mathbb N}$ of diffeomorphisms where the intersections of domains of $F_{i,\delta}$ with $\Pi_\delta$ are pairwise disjoint   and so are the intersections of ranges of $F_{i,\delta}$ with $\Pi_\delta$.


\begin{defi}[Nearly-affine blender system]\label{defi:NABS}
The family $\{\Pi_\delta,\mathcal{F_\delta}\}$ is called a {\em center-stable nearly-affine blender system (cs-NABS)} if the charts $\{(U_\delta,\varphi_\delta)\}$ can be chosen such that, for any $F_{i,\delta}\in \mathcal{F}_\delta$ and  for any point 
$$(x,y,z)\in\psi_\delta(\Pi_\delta)=[-1,1]\times[-1,1]^{d_u}\times[-1,1]^{d_{ss}},\quad 1+d_u+d_{ss}=d,$$
one has that   $(\bar x,\bar y,\bar z)=\varphi_\delta\circ F_{i,\delta}\circ \varphi^{-1}_\delta(x,y,z)$  if and only if 
\begin{equation}\label{eq:sbs:1}
\begin{aligned}
\bar x &= A_{i}(\delta) x + B_{i}(\delta) + \psi_{1,i}(x,\bar{y},z;\delta)=:\hat \psi_{1,i}(x,\bar{y},z;\delta), \\
y &=\psi_{2,i}(x,\bar{y},z;\delta), \\
\bar z &= \psi_{3,i}(x,\bar{y},z;\delta),
\end{aligned}
\end{equation}
where  the coefficients $A_i,B_i$ and functions $\psi$ satisfy that 
\begin{enumerate}[nosep] 
\item[\setword{(S1)}{word:S1}]  $|A_i(\delta)|-\alpha = o(1)_{\delta\to 0}+ o(1)_{\delta \cdot i\to\infty}$ for some constant $\alpha \in (0,1)$. Here it means that given any $\eps>0$, there exists $\delta_0>0$ satisfying that, for any $\delta\in (0,\delta_0)$, one can find  $N(\delta)\in \mathbb N$ such that $|A_i(\delta)|-\alpha<\eps$ for all $i>N(\delta)$,
\item[\setword{(S2)}{word:S2}] for every fixed   $\delta$, the set $\{B_i(\delta)\}_{i\in\mathbb{N}}$ is dense in an interval $I$, containing $0$ and independent of $\delta$ and $i$, and
\item[\setword{(S3)}{word:S3}] the maps $\psi_{j,i}$ $(j=1,2,3)$ are defined on $\{(x,\bar y,z)\in\Pi\}$ and  $\|\psi_{j,i}\|_{C^1}=o(1)_{\delta\to 0}+o(1)_{\delta \cdot i\to\infty}$. 
\end{enumerate}
The family $\{\Pi_\delta,\mathcal{F_\delta}\}$   is called a {\em center-unstable nearly-affine blender system (cu-NABS)} if~\ref{word:S1} is replaced by
\begin{itemize}[nosep] 
\item[\setword{(S1')}{word:S1'}] $|A_i(\delta)|- \alpha' =o(1)_{\delta\to 0}+ o(1)_{\delta \cdot i\to\infty}$ for some constant $\alpha' \in (1,\infty)$.
\end{itemize}
It is called a {\em full nearly-affine blender system} if it is the union of a cs-NABS and a cu-NABS.
\end{defi}


 Let us show that a cs-NABS defines infinitely many cs-Markov partitions. For brevity, we identify $\Pi_\delta$ with $[-1,1]^d=:\Pi$ and     $F_{i,\delta}$ with $\varphi_\delta\circ F_{i,\delta}\circ \varphi^{-1}_\delta$.
 Denote by $\Pi_{i}$ the intersection of $\Pi$ with the domain of $F_{i,\delta}$;  by \eqref{eq:sbs:1} we have
  \begin{equation}\label{eq:domain}
\Pi_{i}=\{(x,y,z): x\in [-1,1], y\in \psi_{2,i}(x, [-1,1]^{d_u},z;\delta),z\in [-1,1]^{d_{ss}}\}.
\end{equation}
  For a finite subset $\mathcal{I} \subset \mathbb{Z}$, define  $g_\delta:\bigcup_{i\in\mathcal{I}}\Pi_{i}\to \Pi$ by $g_\delta(x,y,z)=F_{i,\delta}(x,y,z)$ if $(x,y,z)\in \Pi_{i}$.   Up to shrinking $\Pi$ and rescaling it to the unit size, we can assume that $g_\delta$ is a  diffeomorphism on a small neighborhood of $\bigcup_{i\in\mathcal{I}}\Pi_{i}$. Using \eqref{eq:sbs:1}, \ref{word:S1} and \ref{word:S3}, one readily concludes the following
\begin{prop}
For any sufficiently small $\delta$, there exist  $N_\delta>0$ such that, for any finite subset $I\subset \{i \in \mathbb{Z} : i> N_\delta \}$,  the pair $(\{\Pi_{i}\}_{i\in\mathcal{I}},g_\delta)$ is a cs-Markov partition. We call it a cs-Markov partition of $\{\Pi_\delta,\mathcal{F}_\delta\}$ with specification pair $(\delta,\mathcal{I})$.
\end{prop} 
Example~\ref{ex} illustrates the case where $\card \mathcal{I}=2$. The  cone fields in \eqref{eq:cone} can be obtained by direct computations with formula \eqref{eq:sbs:1}, where the cone coefficients satisfy (see  \citep[Lemma 3.1]{LiTur:24})
\begin{equation}\label{eq:conecoef}
K:=\max\{K^u_i,K^{ss}_i\} =o(1)_{\delta\to 0}+o(1)_{\delta \cdot i\to\infty}.
\end{equation}


When the NABS $\{\Pi_\delta,\mathcal{F_\delta}\}$ is center-unstable, we can also define cu-Markov partitions from the same procedure as above. To see this, we use the Implicit Function Theorem to find $x$ as a function of 
$(\bar x,\bar y, z)$ from the first equation of \eqref{eq:sbs:1}:
$$
x = \dfrac{1}{A_i(\delta)}\bar x - B_i(\delta) + \tilde{\psi}_{1,i}(\bar x,\bar y, z),
$$
for some $\tilde{\psi}_{1,i}$ with $\|\tilde{\psi}_{1,i}\|_{C^1}=o(1)_{\delta\to 0}+o(1)_{\delta \cdot i\to\infty}$. Then, one has that $( x, y, z)=\varphi_\delta\circ F^{-1}_{i,\delta}\circ \varphi^{-1}_\delta(\bar x,\bar y,\bar z)$  if and only if the above equation and the last two of \eqref{eq:sbs:1} are satisfied. Clearly, this defines a cs-NABS $\{\Pi_\delta,\mathcal{F}'_\delta\}$ with $\mathcal{F}'_\delta=\{F^{-1}_{i,\delta}\}$, and cs-Markov partitions of this NABS are cu-Markov partition of the original one.



\begin{thms}\label{thm:nbs}
Let $\{\Pi_\delta,\mathcal{F}_\delta\}$ be a nearly-affine blender system. Up to replacing $\mathcal{F}_\delta$ by some iteration $\mathcal{F}^n_\delta:=\{F^n_{i,\delta}\}$, for any integer $k\geq 1$,  there exists a specification $(\delta,\mathcal{I})$ such that the corresponding partially hyperbolic Markov partition induces a standard blender which is simultaneously separated and $k$-arrayed. It is center-stable/-unstable if   $\{\Pi_\delta,\mathcal{F}_\delta\}$ is center-stable/-unstable.
\end{thms}

We postpone the proof of this theorem to Section~\ref{sec:NBSarray}, after we prove Theorem~\ref{thm:hdc}. 


\subsection{NABS near heterodimensional cycles. Proof of Theorem~\ref{thm:hdc}}\label{sec:LT}

In this section, we prove Theorem~\ref{thm:hdcfull} below, which implies Theorem~\ref{thm:hdc}. Recall that, by saying  a transitive hyperbolic set $\Lambda$ is  $C^1$-wild, we mean that  it exhibits uncountably many $C^1$-robust homoclinic tangencies, and $\homo$ denotes the homoclinic relation between two hyperbolic sets.  

\begin{defi}[$C^1$-robust heterodimensional dynamics/$C^1$-cycle]
We say that a diffeomorphism $f$ has {\em $C^1$-robust heterodimensional cycles} involving two hyperbolic basic sets $\Lambda$ and $\Sigma$ if there exists a $C^1$ neighborhood $\mathcal{U}$ of $f$ such that $W^u(\Lambda_g)\cap W^s(\Sigma_g)\neq \emptyset$ and $W^u(\Sigma_g)\cap W^s(\Lambda_g)\neq \emptyset$ for every $g\in\mathcal{U}$. In this case, we say that  $\Lambda$ and $\Sigma$ form a {\em $C^1$-cycle}.
\end{defi}

Let $f\in \diff^r(\mathcal{M})$, $r=1,\dots,\infty,\omega$,  have a coindex-1 heterodimensional cycle involving two saddles $O_1$ and $O_2$ with $\ind(O_1)+1=\ind(O_2)$.
Recall that we denote by $\lambda$ the center-stable multiplier of $O_1$  and by $\gamma$ the center-unstable multiplier of $O_2$.

\begin{thms}\label{thm:hdcfull}
Given any integer $k> 1$, $f$ can be approximated in the $C^r$ topology by a diffeomorphism $g$, which has up to two  standard blenders, a center-stable one $\Lambda_1$ and a center-unstable one $\Lambda_2$, each of which is simultaneously separated and $k$-arrayed with $\Lambda_i\homo O_{i,g}$ $(i=1,2)$. More specifically,
\begin{itemize}[nosep]
\item in the saddle case, 
\begin{itemize}[nosep]
\item when $|\alpha|\leq 1$, $\Lambda_1$ exists, forming  a $C^1$-cycle with a non-trivial hyperbolic basic set $\Lambda'_2 \homo O_{2,g}$, and
\item when $|\alpha|\geq 1$, $\Lambda_2$ exists, forming  a $C^1$-cycle with a non-trivial hyperbolic basic set $\Lambda'_1 \homo O_{1,g}$;
\end{itemize}
\item in the saddle-focus case, both $\Lambda_1$ and $\Lambda_2$ exist, forming a $C^1$-cycle, and,
\begin{itemize}[nosep]
\item $\Lambda_1$ is $C^1$-wild if  $\gamma$ is nonreal, and
\item $\Lambda_2$ is $C^1$-wild if $\lambda$ is nonreal;
\end{itemize}
\item in the double-focus case, both $\Lambda_1$ and $\Lambda_2$ exist, forming a $C^1$-cycle, and  either $\Lambda_1$  or $\Lambda_2$ is $C^1$-wild.
\end{itemize}
In all cases, $\Lambda_i$ is always $C^1$-wild if $O_i$ originally has a homoclinic tangency.
\end{thms} 

To prove the theorem, we now summarize the  renormalization results in \citep{LiTur:24}, which indicate that each heterodimensional cycle is associated with an NABS.  
Take a {\em fragile} heteroclinic orbit $ \mathcal{T}^0$ from the non-transverse intersection $W^u(O_1)\cap W^s(O_2)$ and a {\em robust} heteroclinic orbit $ \mathcal{T}^1$ from the transverse intersection $W^u(O_2)\cap W^s(O_1)$. We consider the closed invariant set
\begin{equation}\label{cycle}
\Gamma:=\mathcal{O}(O_1)\cup\mathcal{O}(O_2)\cup \mathcal{T}^0 \cup  \mathcal{T}^1,
\end{equation}
where $\mathcal{O}(\cdot)$ denotes the orbit of a point. In \citep{LiTur:24}, a heterodimensional cycle involving $O_1$ and $O_2$ is defined as the set $\Gamma$. To facilitate the presentation, we will do so in the rest of the section.

The key  result in \citep{LiTur:24} regarding blenders is the renormalization of the first return maps for a non-degenerate cycle $\Gamma$ of a $C^r$ ($r=2,\dots,\infty,\omega$) diffeomorphism, where the non-degeneracy conditions  can be fulfilled by an arbitrarily $C^r$-small perturbation. Each first return map is a composition of four maps: two local maps near the  orbits of $O_1$ and $O_2$ which are the restrictions of the period maps, and two transition maps along the heteroclinic orbits $ \mathcal{T}^0$ and $ \mathcal{T}^1$, connecting neighborhoods of two points in $ \mathcal{T}^{0,1}$ with one near $O_1$ and the other near $O_2$.
By considering different iteration time near the two saddle orbits, we obtain a family of first return maps which after renormalization form an NABS. Whether this NABS is center-stable or center-unstable is determined by the  contraction/expansion in the central dynamics. As discussed above Theorem~\ref{thm:hdc}, the central dynamics for saddle cycles is tightly related to the quantity $\alpha$, which in fact appears as
certain first partial derivative of the transition map $\mathcal{T}^1$. In particular, $|\alpha|\neq 1$ for non-degenerate saddle cycles and the resulting NABS   satisfies  $A_i\to \alpha$ in \eqref{eq:sbs:1}.
 We refer the readers to \citep[Sections 2.1--2.3]{LiTur:24} for a detailed description of the first return maps and  non-degeneracy conditions.

\begin{thms}[\citep{LiTur:24}]\label{thm:LT}
Let $f\in \diff^r(\mathcal{M})$, $r=2,\dots,\infty,\omega,$ have a non-degenerate heterodimensional cycle $\Gamma$ of coindex 1, and let  $U$ be any  neighborhood of $\Gamma$. Then, up to an arbitrarily $C^r$-small perturbation  which does not destroy $\Gamma$ when it is of saddle-focus or double-focus type, there exist a family $\{\Pi_\delta\}_\delta$ of embedded cubes in $U$ and  sequences $\{n_i(\delta)\}_{i,\delta}$ of positive integers which are first return times
\footnote{Namely,  $f^n(f^{-n_i(\delta)}(\Pi_\delta)\cap\Pi_\delta)\cap \Pi_\delta=\emptyset$ for $0<n<n_i(\delta)$.}
of $\Pi_\delta$ such that  the family $\{\Pi_\delta,\mathcal{F_\delta}\}$  with  $\mathcal{F_\delta}=\{F_{i,\delta}:=f^{n_i(\delta)}|_{f^{-n_i(\delta)}(\Pi_\delta)}\}$ is a nearly-affine blender system. This system is
\begin{itemize}[nosep]
\item either center-stable or center-unstable if $\Gamma$ is of saddle type, depending on whether $|\alpha|<1$ or $|\alpha|>1$, and
\item full if  $\Gamma$ is of saddle-focus or double-focus type.
\end{itemize}
Any standard cs-blender $\Lambda_1$ and  cu-blender $\Lambda_2$ arisen from the NABS satisfy, up to a further  arbitrarily $C^r$-small perturbation,  that
\begin{itemize}[nosep]
\item in the saddle case,
\begin{itemize}[nosep]
\item if $|\alpha|<1$, either $\Lambda_1\homo O_1$ or it forms a $C^1$-cycle with a non-trivial hyperbolic basic set $\Lambda'_2 \homo O_2$,
\item if $|\alpha|>1$, either $\Lambda_2\homo O_2$ or it forms a $C^1$-cycle with a non-trivial hyperbolic basic set $\Lambda'_1 \homo O_1$, and,
\end{itemize}
moreover, if either $O_1$ or $O_2$ is contained in a non-trivial hyperbolic set, then, in each of the above two cases, the two results hold simultaneously; and
\item in the saddle-focus and double-focus cases,    $\Lambda_1$ and $\Lambda_2$ form a $C^1$-cycle with $\Lambda_1\homo O_1$ and $\Lambda_2\homo O_2$, where the former homoclinic relation can be obtained without destroying the original heterodimensional cycle if $\gamma$ is nonreal, and the same for the latter if $\lambda$ is nonreal.
\end{itemize}
\end{thms}


This theorem follows from several results in \citep{LiTur:24} and is not stated explicitly  there. We detail its references in the \ref{sec:appen}. We also need the following result: 
\begin{lem}\label{lem:HDCtoHT}
Let $O_1$ and $O_2$ have a heterodimensional cycle of coindex-1, and let at least one of the central multipliers $\lambda$ and $\gamma$ be nonreal. Then, there exists an arbitrarily $C^r$-small $(r=1,\dots,\infty,\omega)$ perturbation which yields  a quadratic homoclinic tangency
\begin{itemize}[nosep]
\item to the continuation of $O_1$ if $\gamma$ is not real; and
\item to the continuation of $O_2$ if $\lambda$ is not real.
\end{itemize}
\end{lem}
Note that, even if both $\lambda$ and $\gamma$ are not real, only one of the  two conclusions above can be guaranteed. This lemma can be proved by  computations similar to those  in   \citep[Lemma 2]{GonTurShi:07} and \citep[Lemma 5.3]{LiLiShiTur:22}.
A full proof will appear in the forthcoming paper \citep{BarDiaLi:xx}.
For completeness, we briefly sketch the main idea. We only discuss the first conclusion, as the second one can be achieved by considering $f^{-1}$. By the definition of a heterodimensional cycle, we can find two pieces $W^u\subset W^u(O_1)$ and $W^s\subset W^s(O_1)$ near $O_2$, intersecting $W^s_{loc}(O_2)$ and $W^u_{loc}(O_2)$, respectively. Since $\gamma$ is complex, the forward iterates of $W^u$ are rotating in the central direction, and, in particular, the rotation angles form a dense in $2\pi$ when the argument of $\gamma$, say $\omega$, is irrational. When the cycle is in general position, one can find an iterate of $W^u$ such that it has a tangent vector with arbitrarily small angle to some tangent vector of $W^s$. Thus, a tangency can be created by unfolding the cycle and changing the argument of $\lambda$, say $\omega$, at the same time. Indeed, it is proven in \citep{BarDiaLi:xx} that, if the cycle is non-degenerate, then the family $\{f_{\mu,\omega}\}$, where $\mu$ denotes the signed distance between $W^u$ and $W^s_{loc}(O_2)$ (measured near some point in their original intersection) and $f_{0,\omega^*}=f$, has a sequence $(\mu_i,\omega_i)\to (0,\omega^*)$ of parameter values for which the continuation of $O_1$ has a homoclinic tangency.

\begin{proof}[Proof of Theorem~\ref{thm:hdcfull}]
We prove for the three cases separately.   For brevity, by a perturbation we always mean an arbitrarily $C^r$-small one, and we use the same notations for the continuations of hyperbolic sets after  a perturbation.

\noindent(\textit{Saddle-focus case}.) The existence of the desired  blenders, their homoclinic relations to $O_1$ and $O_2$, and the fact that they form a $C^1$-cycle follow directly from Theorem~\ref{thm:LT}. We thus only prove the part regarding the existence of robust homoclinic tangencies. By considering $f^{-1}$, it suffices to prove that $\Lambda_1$ is $C^1$-wild if either
\begin{enumerate}[nosep]
\item[(1)] $O_1$ has a homoclinic tangency, or
\item[(2)] $\gamma=|\gamma|e^{i\omega}$ with $\omega\in (0,\pi)$.
\end{enumerate}

 Let us start with case (1). By the multidimensional Newhouse theorem, we first, while keep the heterodimensional cycle, unfold this tangency  to obtain a   hyperbolic basic set $\Lambda_0$ which has a $C^2$-robust homoclinic tangency and, by  \citep{GonTurShi:93c,PalVia:94}, is homoclinically related to  $O_1$. Next, we use a further perturbation to make the heterodimensional cycle non-degenerate, and then apply Theorems~\ref{thm:LT} and~\ref{thm:nbs}. After that, we obtain a separated and $k$-arrayed  blender $\Lambda$   which is homoclinically related to $O_1$, and hence also homoclinically related to $\Lambda_0$. By the lambda lemma, the leaves of the stable/unstable invariant  laminations of $\Lambda$ accumulate on those of $\Lambda_0$.  So, one can perturb $f$ again to create a  homoclinic tangency to $\Lambda$, while keeping $\Lambda$ as an separated and $k$-arrayed blender (since this is a $C^1$-open property). Finally, we  invoke Theorem~\ref{thm:uncount} to conclude the proof for the case where $O_1$ has a homoclinic tangency.
 
We proceed to consider case (2), where the center-unstable multiplier $\gamma$ of $O_2$ is nonreal. By Theorems~\ref{thm:LT} and~\ref{thm:nbs}, we first perturb $f$ to obtain a  standard blender which is both separated and $k$-arrayed, and is homoclinically related to $O_1$. Since this perturbation does not destroy the heterodimensional cycle, we apply Lemma~\ref{lem:HDCtoHT} to  create a homoclinic tangency to $O_1$ by a further perturbation. The proof is then completed with  an application of the lambda lemma and  Theorem~\ref{thm:uncount}.

\noindent(\textit{Double-focus case.}) The proof is the same as for the saddle-focus case. 

\noindent(\textit{Saddle case}.) The proof in this case is almost identical to case (1) of the saddle-focus case, except that the $C^1$-cycle involving $\Lambda_1$ and $\Lambda'_2$ is not automatic. To obtain this, we need to use the `moreover' part of the saddle case in Theorem~\ref{thm:LT}, namely, we need to show that $O_1$ or $O_2$ is contained in non-trivial hyperbolic set. But this is immediate, because, after the first step of case (1), either $O_1$ or $O_2$ (depending on which one has the tangency) is homoclinically related to the hyperbolic basic set $\Lambda_0$.
\end{proof}

\subsection{Standard blenders arisen from an NABS. Proof of Theorem~\ref{thm:nbs}}\label{sec:NBSarray}

We will only consider the  cs-NABS, and the proof for the center-unstable case is completely parallel.  Let us identify  $F_{i,\delta}$ with $\varphi_\delta\circ F_{i,\delta}\circ \varphi^{-1}_\delta$. So, $F_{i,\delta}$ is a map defined on  $\Pi_i$ given by \eqref{eq:domain} and satisfy the cross-form relation \eqref{eq:sbs:1}. 
We prove Theorem~\ref{thm:nbs} in several steps. Roughly speaking, we first find  a specification set $\mathcal{I}$ such that   the BCG covering property in Definition~\ref{defi:bcgcover} is satisfied, and hence the cs-Markov partition $(\{\Pi_{i}\}_{i\in\mathcal{I}},g_\delta)$ induces is a $k$-arrayed  blender. Then, we show that by `duplicating' this set $\mathcal{I}$ the obtained arrayed blender can also be a separated one.

\subsubsection{Preliminaries.}  Let us define the BCG structure  for $\Pi_i$ in Definition~\ref{defi:structure2} (which also includes the BC structure in Definition~\ref{defi:BC}). By \eqref{eq:domain}, the sets $\mathbb{X}_i,\mathbb{Y}_i,\mathbb{Z}_i$ in Definition~\ref{defi:markov} are the same for all $i$ in our case, and, after dropping the subscript, they are $\mathbb{X}=[-1,1],\mathbb{Y}=[-1,1]^{d_u},\mathbb{Z}=[-1,1]^{d_{ss}}$. Denote $\Pi=[-1,1]^d$.
 We take real numbers $x_\B,x_\C,x_\G$ with $0<x_\C<x_\G<x_\B<1$ to be determined, and define  for all $i$
$$\mathbb{X}^\B=(-x_\B,x_\B),\quad \mathbb{X}^\C=(-x_\C,x_\C),\quad
x^\L=-x_\G,\quad x^\R=x_\G,$$
and
\begin{equation}\label{eq:bcg}
\begin{aligned}
&\Pi^\B=(-x_\B,x_\B)\times \mathbb{Y} \times\mathbb{Z},\quad
\Pi^\C=(-x_\C,x_\C)\times \mathbb{Y} \times\mathbb{Z},\\
&\mathring\Pi^\L = \{-x_\G\}\times \mathbb{Y} \times\mathbb{Z},\quad
\mathring\Pi^\R = \{x_\G\}\times \mathbb{Y} \times\mathbb{Z},\\
&\Pi_i^\B=\Pi_i\cap \Pi^\B,\quad
\Pi_i^\C=\Pi_i\cap \Pi^\C,\quad
\mathring\Pi^\L_i =  \Pi_i\cap \mathring\Pi^\L,\quad
\mathring\Pi^\R_i =  \Pi_i\cap \mathring\Pi^\R.
\end{aligned}
\end{equation}

 Consider the map $g$ defined on all $\Pi_i$ by
$g(x,y,z)=F_{i,\delta}(x,y,z)$ if $(x,y,z)\in \Pi_i$. By \eqref{eq:sbs:1}, we have
 \begin{equation}\label{eq:V}
 \begin{aligned}
V_i&:=g(\Pi_i)= \{(\hat\psi_{1,i}(x,y,z),y,\psi_{3,i}(x,y,z)):(x,y,z)\in \X\times\Y\times \Z\},\\
V_i^{\B/\C}&:=g(\Pi^{\B/\C}_i)= \{(\hat\psi_{1,i}(x,y,z),y,\psi_{3,i}(x,y,z)):(x,y,z)\in \X^{\B/\C}\times\Y\times \Z\},\\
\mathring{V}_i^{\L/\R}&:=g(\mathring{\Pi}_i^{\L/\R})= \{(\hat\psi_{1,i}(x,y,z),y,\psi_{3,i}(x,y,z)):(x,y,z)\in \{x^{\L/\R}\}\times\Y\times \Z\}.
 \end{aligned}
\end{equation}
Note that  the vertical strips which constitute the vertical arrays in the BCG covering property (Definition~\ref{defi:bcgcover}) are those inside $\Pi_j$. Hence, in our case, they are of the form $V_{ij}=V_i\cap \Pi_j$. Similarly, we have $V^{\B/\C/\L/\R}_{ij}=V^{\B/\C/\L/\R}_i\cap \Pi_j$. Combining this fact with the above BCG structure of $\Pi_i$, we obtain
\begin{lem}\label{lem:equiv}
The  hyperbolic basic set $\Lambda$ induced from the cs-Markov partition $(\{\Pi_{i}\}_{i\in\mathcal{I}},g_\delta)$ is a $k$-arrayed  blender if  the BCG covering property is satisfied with $\Pi$ and $\mathcal{V}\subset\{V_i:i\in\mathcal{I}\}$, in place of $\Pi_i$ and $\mathcal{V}\subset\{V_{ij}:i,j\in\mathcal{I}\}$.
\end{lem}

 Recall Definition~\ref{defi:disc} that an ss-disc $S$ is of the form $(x,y)=(s_1(z),s_2(z))$ for $z\in \mathbb{Z}$ and a u-disc $S^u$ is of the form $(x,z)=(s^u_1(y),s^u_2(y))$ for $y\in \mathbb{Y}$. Let $\Delta^{ss}_S$ and $\Delta^{u}_{S^u}$ be the maximal deviation of the $x$-coordinate of an ss-disc $S$ and a u-disc $S^u$, i.e., 
$$\Delta^{ss}_S =\max_{z_1,z_2\in[-1,1]^{d_{ss}}}|s_1(z_1)-s_1(z_2)|
\quad\mbox{and}\quad
\Delta^{u}_{S^u} =\max_{y_1,y_2\in[-1,1]^{d_{ss}}}|s^u_1(y_1)-s^u_1(y_2)|.
$$
Denote by $\Delta_i^{ss}$  the supremum of $\Delta^{ss}_S $  over all ss-discs $S$   crossing $\Pi_i$, and by $\Delta_i^u$ the supremum of $\Delta^{u}_{S^u}$ over all u-discs $S^u$ crossing $\Pi_i$.  
One sees from \eqref{eq:conecoef} that $\max\{\Delta_i^{ss},\Delta_i^u\}=o(1)_{\delta}+o(1)_{\delta\cdot i \to\infty}$. (Recall the meaning of this notation in \ref{word:S1}.)

By  \eqref{eq:sbs:1}, the image $\hat\psi_{1,i}(\X^{\B/\C}\times\Y\times \Z)$ contains the interval $(B_i-A_ix_\B+\| \psi_{1,i}\|,B_i+A_ix_\B-\| \psi_{1,i}\|)$. By the form of $V^\B_i$ in \eqref{eq:V}, one can easily check\footnote{An elaborated computation can be found in the proof of \citep[Lemma 4.10]{LiTur:24}.} that  any ss-disc which lies entirely in 
\begin{equation}\label{eq:cross2}
(B_i-|A_i|x_\B+\| \psi_{1,i}\|,B_i+|A_i|x_\B-\| \psi_{1,i}\|)\times \Y\times \Z
\end{equation}
 crosses $V^\B_i$ in the sense of Definition~\ref{defi:base}. By \ref{word:S1}, \ref{word:S3} and \eqref{eq:conecoef}, we have
 \begin{equation}\label{eq:delta} 
 \Delta_i:=\max\{ |A_i|-\alpha, \| \psi_{1,i}\|,\Delta_i^{ss},\Delta_i^u\}
 =o(1)_{\delta}+o(1)_{\delta\cdot i \to\infty}.
 \end{equation}

Denote
\begin{equation}\label{eq:XV^B}
\begin{aligned}
\mathbb{X}(V^\B_i)=(B_i-\alpha x_\B+3\Delta_i,B_i+\alpha x_\B-3\Delta_i),\\
\mathbb{X}(V^\C_i)=(B_i-\alpha x_\C+3\Delta_i,B_i+\alpha x_\C-3\Delta_i).
\end{aligned}
\end{equation}
\begin{lem}\label{lem:cross}
If an ss-disc  crosses $\Pi^\B$ and intersects $\mathbb{X}(V^\B_i)\times \Y \times \Z$, then it crosses $V^\B_i$. Similarly, if it intersects $\mathbb{X}(V^\C_i)\times \Y \times \Z$, then it crosses $V^\C_i$.
\end{lem}
\begin{proof}
By Definition~\ref{defi:base} of crossing, the ss-disc lies entirely in $(-x_\B,x_\B)\times \Y \times \Z$. By the choice of $\Delta$, if it intersects the cube in the statement of the lemma, then it must lie entirely in the cube given by \eqref{eq:cross2}.
\end{proof}

\subsubsection{Finding 1-arrayed standard blenders.}
Let us now find a specification set $\mathcal{I}$  such that $(\{\Pi_{i}\}_{i\in\mathcal{I}},g_\delta)$ satisfies 1-arrayed BCG covering property as in Definition~\ref{defi:bcgcover}, i.e., every vertical array  has only one vertical strip. After that, we show that the $k$-arrayed covering property can be satisfied by repeating the same argument $k$ times. By Lemma~\ref{lem:equiv}, we verify the properties~\ref{word:B1}--\ref{word:B4} with $\Pi$ and $\mathcal{V}\subset \{V_i:i\in\mathcal{I}\}$ in Lemmas~\ref{lem:bcg1}--\ref{lem:bcg5} below.  


Take any $x_\B\in(0,1)$ such that $[-x_\B,x_\B]$ is contained in the interval in which $\{B_i\}$ is dense, given by \ref{word:S2}.
Let $N>0$ be an integer to be determined.  We divide the interval $[-x_\B,x_\B]$  equally into $N$ subintervals and denote the endpoints by $x_n$ with $n=0,\dots,N$, $x_0=-x_\B$ and  $x_N=x_\B$.  By \eqref{eq:delta}, we can take a sequence $\{\delta_s\}\to 0$ and a family  $\{\mathcal{I}_s\}_{s\in\mathbb{N}}$ of sets $\mathcal{I}_s=\{i_{1}(s),\dots,i_{N}(s)\}$ such that
$$
\Delta_{i_{n}(s)}\to 0
\spc
B_{i_{n}(s)}-x_n\to 0
\quad\mbox{as}\quad
s\to\infty.
$$
We further require the family of pairs $\{(\delta_s,\mathcal{I}_s)\}$ is chosen such that $\Delta_{i_n(s)}\to 0$ as $s\to\infty$.
Let $\floor{\cdot}$ denote the largest integer smaller than a given number.

\begin{lem}\label{lem:bcg1}
Given any  integer $N$ with $N\geq \floor{1/\alpha}+2$, the specification sets $\mathcal{I}_s$ satisfy  property~\ref{word:B1} for all  sufficiently large $s$.
\end{lem}

\begin{proof}
Note that, if the union of the associated intervals $\X(V^\B_i)$ with $i\in\mathcal{I}$ covers $[-x_\B,x_\B]$, then every ss-disc  crossing $\Pi^\B$ must intersect $\mathbb{X}(V^\B_{i})\times \Y\times \Z$ for some $i$, and hence cross $V^\B_{i}$ by Lemma~\ref{lem:cross}, which gives~\ref{word:B1}. 

Taking
\begin{equation}\label{eq:deltaa}
\tilde{\Delta}_s=\max_{n=1,\dots,N} \{\Delta_{i_n(s)},|B_{i_{n}(s)}-x_n|\}, 
\end{equation}
one sees that   $\X(V^\B_{i_{n}(s)})$  contains the subinterval
$$
\tilde \X(V^\B_{i_{n}(s)}):= (x_n-\alpha x_\B+4\tilde\Delta_s,x_n+\alpha x_\B-4\tilde\Delta_s).
$$
Hence, $[-x_\B,x_\B]$ is covered by $\bigcup_n \X(V^\B_{i_{n}(s)})$ if it is covered by $\bigcup_n\tilde \X(V^\B_{i_{n}(s)})$. But the latter is obvious for all sufficiently large $s$, since 
$$
\sum_n \left|\tilde \X(V^\B_{i_{n}(s)})\cap[-x_\B,x_\B]\right|
>(2\alpha x_\B-8\tilde\Delta_s)\cdot (\alpha^{-1}+1)=2 x_\B+2\alpha x_\B+o(1)_{ s\to\infty},
$$
where $|\cdot |$ denotes the length of an interval.  It then follows from the discussion at the beginning that~\ref{word:B3} holds for all $\mathcal{I}_s$ with  sufficiently large $s$.
\end{proof}

\begin{lem}\label{lem:bcg2}
The specification sets $\mathcal{I}_s$ satisfy property~\ref{word:B4}  for all  sufficiently  large $s$.
\end{lem}
\begin{proof}
By the choice of $\Delta_i$ in \eqref{eq:delta}, Definition~\ref{defi:structure2} and the BCG structure in \eqref{eq:bcg}, we find that
\begin{equation}\label{eq:PiLR}
\Pi^{\L}\subset (-x_\G-3\Delta_i,-x_\G+3\Delta_i)\times \Y\times\Z
\spc
\Pi^{\R}\subset (x_\G-3\Delta_i,x_\G+3\Delta_i)\times \Y\times\Z,
\end{equation}
where we use three $\Delta_i$  to bound $\|\psi_{1,i}\|$, the $x$-deviations of the ss-discs and of the u-discs, respectively. Recall the separation notation $\parallel$ is defined before Definition~\ref{defi:bcgcover}. It is immediate from \eqref{eq:bcg} that  $\Pi^{\L/\R}_i\parallel \Pi_i^{\C}$ if $x_\G-x_\C>4\Delta_i$, which is automatic for all sufficiently large $i$. 

Let us check $V^{\R}_i\parallel V^{\C}$. Recall that $|A_i-\alpha|<\Delta_i$.  We find from \eqref{eq:sbs:1} and \eqref{eq:V}  that 
\begin{align*}
&V_i^{\R}
\subset ((\alpha+\Delta_i)(x_\G-3\Delta_i)+B_i-\Delta_i,(\alpha+\Delta_i)(x_\G+3\Delta_i)+B_i+\Delta_i)\times \Y\times\Z,\\
&V_i^\C\subset
(-(\alpha+\Delta_i)x_\C+B_i-\Delta_i,(\alpha+\Delta_i)x_\C+B_i+\Delta_i)\times \Y\times\Z.
\end{align*}
Thus, $V^{\R}_i\parallel V^{\C}$ if 
$$
((\alpha+\Delta_i)(x_\G-3\Delta_i)+B_i-\Delta_i) - ((\alpha+\Delta_i) x_\C+B_i+\Delta_i)>\Delta_i,
$$
where the $\Delta_i$ on the right-hand side is used to bound the $x$-deviation  of the ss-discs. This inequality  is satisfied when   $ x_\G- x_\C> 9\Delta_i \alpha^{-1}$. By a similar discussion, the same condition can be derived from $V^{\L}_i\parallel V^{\C}$ and $V^{\L/\R}_i\parallel \partial^{\L/\R}_x V^{\B}$. The lemma follows immediately since  by construction $i\in \mathcal{I}_s$ tend to infinity as $s\to \infty$.
 \end{proof}

\begin{lem}\label{lem:bcg3}
Properties~\ref{word:B3} and~\ref{word:B3'} are satisfied by $\mathcal{I}_s=\{x_0,\dots,x_N\}$ for all  sufficiently  large $s$ if 
\begin{itemize}[nosep]
\item $N$ is odd and $-x_\C=x_n$ for some $0<n<(N-1)/2$, and
\item $x_\G>x_\C+\alpha x_\B$.
\end{itemize} 
\end{lem}

\begin{proof}
Let  $\tilde{\Delta}_s$ be given by \eqref{eq:deltaa} and recall $\mathbb{X}(V^\C_i)$ defined in \eqref{eq:XV^B}. By \eqref{eq:sbs:1} and \eqref{eq:V}, 
\begin{equation}\label{eq:VC}
\mathbb{X}(V^\C_{i_n(s)})\times \Y\times\Z
\supset (x_n- \alpha x_\C +3\tilde{\Delta}_s,x_n+ \alpha x_\C -3\tilde{\Delta}_s)\times \Y\times\Z,
\end{equation}
where the term  $3\tilde{\Delta}_s$ comes from $|\alpha-A_{i_n(s)}|,|x_\C-B_{i_n(s)}|$, and $\|\psi_{1,i_n(s)}\|$. 
Since every ss-disc intersecting $\partial^\L_x\Pi^\C=\{-x_\C\}\times \Y\times\Z$ has its $x$-deviation bounded by $2\tilde{\Delta}_s$ from $x=-x_\C=x_n$, it must cross the second cube in \eqref{eq:VC}.
$
2\alpha x_\C - 6\tilde{\Delta}_s>2\tilde{\Delta}_s$,
 which is automatic for all large $s$. This together with Lemma~\ref{lem:cross} proves the first part of~\ref{word:B3}.

To check $V^\B_{i_n(s)}\parallel \Pi^\L$, we obtain from \eqref{eq:sbs:1} and \eqref{eq:V} that 
\begin{equation}\label{eq:VB}
V^\B_{i_n(s)}\subset (x_n- \alpha x_\B -3\tilde{\Delta}_s,x_n+ \alpha x_\B +3\tilde{\Delta}_s)\times \Y\times\Z.
\end{equation}
Since $ \Pi^\L$ belongs to the first cube in \eqref{eq:PiLR} and  $x_n=-x_\C$, it follows that $V^\B_{i_n(s)}\parallel \Pi^\L$ if $-x_\C- \alpha x_\B + x_\G> 7\tilde{\Delta}_s,$ which is fulfilled by all large $s$ due to the second condition in the lemma. Since $N$ is odd, the discussion for~\ref{word:B3'} is the same by symmetry.
\end{proof}

Recall that the points $x_n$ divide $[-x_\B,x_\B]$ into $N$ subintervals of the same length.
\begin{lem}\label{lem:bcg4}
Properties~\ref{word:B2} and~\ref{word:B2'} are satisfied by $\mathcal{I}_s=\{x_0,\dots,x_N\}$ for all  sufficiently  large $s$ if 
\begin{equation}\label{eq:B3}
\alpha(x_\B-x_\C)<x_{n+1}-x_n<\alpha(x_\B+x_\C).
\end{equation}
More specifically, for $n=0,\dots,N-1$, any ss-disc intersecting $\partial^{\mathrm{R}}_x V^\B_{i_n(s)}$ crosses $V^\C_{i_{n+1}(s)}$,  and $ V^\B_{i_{n+1}(s)} \parallel V_{i_n(s)}^\L$; and, for $n=1,\dots,N$, any ss-disc intersecting $\partial^{\mathrm{L}}_x V^\B_{i_{n+1}(s)}$ crosses $V^\C_{i_{n}(s)}$,  and $ V^\B_{i_{n}(s)} \parallel V_{i_{n+1}(s)}^\R$.
\end{lem}
\begin{proof}
Again by \eqref{eq:sbs:1} and \eqref{eq:V}, we find that $\partial^{\R}_x V^\B_{i_n(s)}$ lies between the following two hyperplanes:
$$
\{x_n + \alpha x_\B-3\tilde{\Delta}_s\}\times \Y\times\Z
\spc
\{x_n + \alpha x_\B+3\tilde{\Delta}_s\}\times \Y\times\Z,
$$
while, by \eqref{eq:VC}, 
$$
\mathbb{X}(V^\C_{i_{n+1}(s)})\times \Y\times\Z
\supset 
(x_{n+1}- \alpha x_\C +3\tilde{\Delta}_s,x_{n+1}+ \alpha x_\C -3\tilde{\Delta}_s)\times \Y\times\Z,
$$
where $\mathbb{X}(V^\C_{i_{n+1}(s)})$ is defined in \eqref{eq:XV^B}.
Thus, by Lemma \ref{lem:cross}, an ss-disc that intersects $\partial^{\R}_x V^\B_{i_n(s)}$ will cross $V^\C_{i_{n+1}(s)}$ if
\begin{align*}
x_n + \alpha x_\B-3\tilde{\Delta}_s - (x_{n+1}- \alpha x_\C +3\tilde{\Delta}_s)>\tilde{\Delta}_s,\\
x_{n+1}+ \alpha x_\C -3\tilde{\Delta}_s-(x_n + \alpha x_\B+3\tilde{\Delta}_s)>\tilde{\Delta}_s,
\end{align*}
or
\begin{equation*}
x_{n+1}-x_n<\alpha(x_\B+x_\C) - 7\tilde{\Delta}_s
\spc
x_{n+1}-x_n>\alpha(x_\B-x_\C) +  7\tilde{\Delta}_s,
\end{equation*}
which holds for all sufficiently large $s$ by \eqref{eq:B3}. This proves the first part of~\ref{word:B2}. 

To check the second part, we find from \eqref{eq:VB}
$$
V^\B_{i_{n+1}(s)}\subset (x_{n+1}- \alpha x_\B -3\tilde{\Delta}_s,x_{n+1}+ \alpha x_\B +3\tilde{\Delta}_s)\times \Y\times\Z,
$$
and from \eqref{eq:sbs:1} and \eqref{eq:V}
$$V_{i_{n}(s)}^{\L}
\subset (-\alpha x_\G+x_n-7\tilde\Delta_s,-\alpha x_\G+x_n+7\tilde\Delta_s)\times \Y\times\Z.$$
Obviously, we have $ V^\B_{i_{n+1}(s)} \parallel V_{i_n(s)}^\L$ if
$$
x_{n+1}- \alpha x_\B -3\tilde{\Delta}_s
-(-\alpha x_\G+x_n+6\tilde\Delta_s)>\tilde{\Delta}_s,
$$
which holds for all sufficiently large $s$  by \eqref{eq:B3}.
The proof of~\ref{word:B2'} is completely parallel due to the symmetry.
\end{proof}

\begin{lem}\label{lem:bcg5}
If $\alpha < 1/20$, then $x_\C,x_\G,N$ can be chosen such that the conditions in Lemmas~\ref{lem:bcg1},~\ref{lem:bcg3} and \ref{lem:bcg4} are satisfied.
\end{lem}

\begin{proof}
Since $x_{n+1}-x_n=2x_\B/N$, inequality \eqref{eq:B3} can be rewritten as 
\begin{equation}\label{eq:ineN}
\dfrac{2x_\B}{\alpha(x_\B+x_\C)}< N< \dfrac{2x_\B}{\alpha(x_\B-x_\C)}.
\end{equation}
We claim that, for every $N$ satisfying this inequality and has the form $N=8m+1$ with $m\geq 1$, one can define the desired  $x_\C,x_\G$. 

Let us prove the claim.
 Take $x_\C=-x_{m+1}$, so
${x_\C}/{x_\B}=(m-1/2)/(4m+1/2)\in  (1/9,1/4)$, and hence
\begin{equation}\label{eq:N1}
\dfrac{1}{9}x_\B<x_\C<\dfrac{1}{4}x_\B.
\end{equation}
This along with \eqref{eq:ineN} implies
\begin{equation}\label{eq:N2}
\dfrac{8}{5\alpha}<N<\dfrac{9}{4\alpha}.
\end{equation}
The condition of Lemma~\ref{lem:bcg1} is satisfied since
${8}/{(5\alpha)}-(\floor{{1}/{\alpha}}+2)
>32 - 22>0 $. The first condition of Lemma~\ref{lem:bcg3} is automatic by construction, and the second one is equivalent to require $x_\C+\alpha x_\B<1$, which is also satisfied since $x_\C+\alpha x_\B<x_\B/4 + x_\B/20 <1$. So, $x_\C$ can be taken as any number satisfying the second condition of  Lemma~\ref{lem:bcg3}. The claim is proven.

It remains to show that there indeed exist  integers of the form $N=8m+1$ satisfying \eqref{eq:ineN}. But this is obvious since
$$
 \dfrac{2x_\B}{\alpha(x_\B-x_\C)} - \dfrac{2x_\B}{\alpha(x_\B+x_\C)}>
 \dfrac{9}{4\alpha}-\dfrac{9}{5\alpha}>9.
$$
\end{proof}

%
%
%
%

\subsubsection{Finding $k$-arrayed standard blenders.}\label{sec:proofk}
 Lemmas~\ref{lem:bcg1}--\ref{lem:bcg5} imply that, when $\alpha<1/20$, there exists a family $\{\mathcal{I}_s\}$ of specification sets defined above Lemma~\ref{lem:bcg1} such that for  every sufficiently large $s$, the cs-Markov partition $(\{\Pi_i\}_{i\in\mathcal{I}_s},g_\delta)$ induces a 1-arrayed  blender. To obtain a $k$-arrayed  blender for any $k>1$, we find $k$ families $\{\mathcal{I}^j_s\}_{s\in\mathbb{N},1\leq j \leq k}$ via Lemmas~\ref{lem:bcg1}--\ref{lem:bcg5} by taking $\tilde\Delta_s=\max_j\{\tilde\Delta^j_s\}$, where $\tilde\Delta^j_s$ are defined as \eqref{eq:deltaa} but now with the superscript. Then, for every small $\delta$ and large $s$, the cs-Markov partition $(\{\Pi_i\}_{i\in\bigcup_j\mathcal{I}^j_s},g_\delta)$ induces a $k$-arrayed  blender. To remove the assumption $\alpha<1/20$, one simply considers the iteration $g^n$ for some $n$ such that the new contraction coefficient $\alpha^n<1/20$.

\subsubsection{Building a separated blender: completion of the proof.}
Recall that the BC structure (see Definition~\ref{defi:BC}) is defined as in \eqref{eq:bcg}. Denote by $\mathcal{I}$ the specification set of the found $k$-arrayed  blender. In what follows, we show that by adding two more copies of $\mathcal{I}$ with  some shifts in the $x$-direction, the resulting bigger specification set yields a standard blender which is simultaneously $k$-arrayed and separated.

Recall the numbers $x_i$ $(i=0,\dots,N)$ defined above Lemma~\ref{lem:bcg1} and denote $D=x_1-x_0$. We take
$$
x'_i=x_i+\dfrac{1}{3}D, \quad
x_{-1}' = x_0 - \dfrac{2}{3}D,\quad
x''_i=x_i+\dfrac{2}{3}D, \quad
x_{-1}'' = x_0 - \dfrac{1}{3}D.
$$
Let $\mathcal{I}'$ and $\mathcal{I}''$ be the specification sets found in the same way as $\mathcal{I}$, but with $\{x'_i\}$ and $\{x''_i\}$ instead  of $\{x_i\}$  in Lemmas~\ref{lem:bcg1}--\ref{lem:bcg5}. Evidently, the standard blender of $\hat{\mathcal{I}}=\mathcal{I}\cup\mathcal{I}'\cup \mathcal{I}''$ is still a $k$-arrayed one, where the extra points $x'_{-1}$ and $x_{-1}''$ are for ensuring that property~\ref{word:B2'} is satisfied by each of the vertical array sets corresponding to $\mathcal{I}'$ and $\mathcal{I}''$.

We now prove that the blender of $\hat{\mathcal{I}}$ is also a separated one. Property~\ref{word:A1} is immediate from~\ref{word:B1}. To verify~\ref{word:A2}, it suffices to show that, instead of $\mathcal{V}^\B_2\parallel \mathcal V^\L_1 $ in~\ref{word:B2}, the stronger property $\mathcal{V}^\C_2\parallel \mathcal V^\C_1 $ holds, and Similarly for~\ref{word:B2'}. Due to the symmetry, we only prove $\mathcal{V}^\C_2\parallel \mathcal V^\C_1 $. One just needs to replace $V^\B_{i_{n+1}(s)}$ and $V^\L_{i_{n}(s)}$ in the second part of the proof  of Lemma~\ref{lem:bcg4} by $V^\C_{i_{n+1}(s)}$ and $V^\C_{i_{n}(s)}$. This yields a condition that can be satisfied for all large $s$ if $x_{n+1}-x_n>2\alpha x_\C  $, which, by the choice of $x_\C$ in Lemma~\ref{lem:bcg5} and the fact that $x_{n+1}-x_n=2x_\B/N$, is equivalent to 
$$
\alpha <\dfrac{4m+\frac{1}{2}}{(8m+1)(m-\frac{1}{2})}.
$$
This can be satisfied for any fixed $m$ by consider iterations of the NABS, as explained in Section~\ref{sec:proofk}. Thus,~\ref{word:A2} is satisfied.

We proceed to check a stronger version of~\ref{word:A3} that every ss-disc intersecting $\Pi_i^\B$ will cross some $V^\C$. This is just a repetition of the proof of Lemma~\ref{lem:bcg1}, with requiring that $[-x_\B,x_\B]$ is covered by $\mathbb{X}(V)^\C_i$, instead of $\mathbb{X}(V)^\B_i$,  with $i\in \hat{\mathcal{I}}$. It follows that the desired property holds if 
$
\card{\hat{\mathcal{I}}}=3N \geq \floor{{x_\B}/{\alpha x_\C}}+2
$, which is automatic by \eqref{eq:N1} and \eqref{eq:N2}. Thus, the $k$-arrayed  blender of ${\hat{\mathcal{I}}}$ is also separated, and the proof of Theorem~\ref{thm:nbs} is completed.

\appendix
\section{NABS obtained in \citep{LiTur:24}}
\makeatletter\def\@currentlabel{Appendix}\makeatother
\label{sec:appen}

\subsection{Proof of Theorem~\ref{thm:LT}}\label{sec:appen1}
We now collect several results in \citep{LiTur:24} which together imply Theorem~\ref{thm:LT}. In some cases, we will also refer to \citep{Li:24} where the corresponding improved results are obtained. The statements there hold for all standard blenders arisen from the NABS near heterodimensional cycles,  which in particular include the separated and arrayed ones introduced in this paper.

  In \cite[Section 2.5]{LiTur:24}, saddle heterodimensional cycles are further divided into two subcases, type I and type II, according to whether  the transition map along the non-transverse heteroclinic orbit $\mathcal{T}^0$ satisfies an orientation-related condition  or not. The existence of the NABS for type-I  cycles follows from the renormalization formulas (3.26) and (4.21) of \citep{LiTur:24}, with replacing $(X,Y,Z)$ by $(\delta x,\delta y,\delta z)$ in these two formulas.
The part on the homoclinic relation and $C^1$-cycle follows from \citep[Theorem 2]{LiTur:24}, while the additional `moreover' part (concerning the case where one of the saddles is contained in a non-trivial hyperbolic set) is by \citep[Corollary 3]{LiTur:24}.
The NABS part for type-II  cycles is due to  \cite[Theorem 5]{LiTur:24} that  saddle cycles of type II  are accumulated by those of type I,  and other parts are from \citep[Corollary 2]{LiTur:24}. 
The proofs are based on the analysis of first return maps on a family $\{U_\delta\}$ of shrinking neighborhoods of some point $P\in \Gamma$, where $\delta$ denotes the size. A single renormalization map $\varphi:U\to \mathbb{R}^d$  for some $U\supset U_\delta$ is considered, that is, $\varphi_{\delta}=\varphi|_{U_{\delta}}$. In this case, one has $\{n_i(\delta_2)\}\subset \{n_i(\delta_1)\}$ for any $\delta_2<\delta_1$.

 
The NABS part of Theorem~\ref{thm:LT} for  saddle-focus and double-focus cycles follows from formulas (6.18) and, respectively, (6.44) of \citep{LiTur:24}, with  replacing $(X,Y,Z)$ and, respectively, $(U,V,W)$ by $(\delta x,\delta y,\delta z)$. A more concise presentation of these two formulas can  be found in \citep[Sections 2.2 and 2.3]{Li:24}. The convergence requirements in \ref{word:S1} and \ref{word:S1'} are described in detail by Lemma~\ref{lem:fullcon} below.

The part on the homoclinic relation and $C^1$-cycle is given by \citep[Theorem 7]{LiTur:24}. In these two cases, the proofs are still based on a family of shrinking neighborhoods but now the involved renormalization maps depend on the size of the neighborhoods.

\begin{rem}
The first perturbation  in Theorem~\ref{thm:LT} is to achieve certain arithmetic properties on the central multipliers, which guarantee the existence of the nearly-affine blender systems. These properties are fully classified in \citep[Corollary A]{Li:24}.
 The exception is for saddle cycles of type II, where the cycle is additionally unfolded to create those of type I.
 The second perturbation for homoclinic relations and $C^1$-cycles are done by unfolding the  heterodimensional cycle, and this can be done within generic finite-parameter families, where a complete discussion is given in by \citep[Theorem 2]{LiTur:24} for the saddle case and \citep[Theorem C]{Li:24} for other cases.
\end{rem}

\subsection{NABS in the saddle-focus and double-focus cases}
Theorem~\ref{thm:LT} in these two cases establishes the coexistence of a cs-NABS and a cu-NABS, and they carry additional properties   not used in this paper.  First, in contrast to the saddle case, where $A_i$ in \eqref{eq:sbs:1} converge to the quantity $\alpha$, we have
\begin{lem}[Full central contraction/expansion]\label{lem:fullcon}
In the saddle-focus and double-focus cases, it holds that for any $\alpha\in\mathbb R\setminus\{0,\pm 1,\pm\infty\}$, there exists a sequence $\{i_n\}_{n\in\mathbb{N}}$ such that $A_{i_n}\to \alpha$ as $n\to \infty$ and $\{B_{i_n}\}$ is dense near $0$.
\end{lem}

 Recall that the map \eqref{eq:sbs:1} are obtained as the first return map along the heterodimensional cycle, and they are of the form $F_i:=F_{k_i,m_i}$, where $k_i$ and $m_i$ are the  number of iterations spent near $O_1$ and $O_2$, respectively. The result  in the saddle-focus case follows from the proof of \citep[Proposition 6.2]{LiTur:24}, or, more specifically, the last displayed equation in the proof. For convenience we present it here with a concise form.
$$A_{k_i,m_i}=\dfrac{(B_{k_i,m_i}-c)\sin(k\omega+\eta_1)}{\sin(k\omega+\eta_2)},$$
where $A_{k_i,m_i},(B_{k_i,m_i}-c)$ are the quantities $A_i$ and $B_i$ in \eqref{eq:sbs:1}, $\omega$ is the argument of the complex center-stable multiplier of $O_1$, and $c,\eta_{1,2}$ are constants with $\eta_1\neq \eta_2+n\pi$. The arithmetic condition of \citep[Theorem 7]{LiTur:24}  ensures that $\{k\omega/2\pi, B_{k,m}\}_{k,m\in\mathbb{N}}$ is dense in $\mathbb{R}^2$. This implies that, given any $\alpha$, one can choose  $(k_i,m_i)$ such that $\{(B_{k_i,m_i}-c)\}$ is dense near a small neighborhood of $0$, and, in the same time, $A_{k_i,m_i}\to \alpha$, proving the lemma in the case where the center-stable multipliers of $O_1$ are nonreal. When $O_2$ has nonreal center-unstable multipliers, one just considers $f^{-1}$. This lemma for the double-focus case follows from a reduction to the saddle-focus one in the proof of  \citep[Proposition 6.7]{LiTur:24}.

\begin{lem}[Overlapping]
The cs- and cu-NABSs, which constitute the full NABS in the saddle-focus and double-focus cases, share the same family $\{(U_\delta,\pp_\delta)\}$ of charts in Definition~\ref{defi:NABS}.
\end{lem}

This lemma follows directly from \citep[Propositions 6.2 and 6.7]{LiTur:24}. An immediate consequence due to Proposition~\ref{prop:nontr} is that, for any pair of standard cs- and cu-blender arisen from such full NABS, denoted as $\Lambda^{cs}$ and $\Lambda^{cu}$, one has that $W^u(\Lambda^{cs})$ intersects $W^s(\Lambda^{cu})$ $C^1$-robustly.

\section*{Acknowledgment}
I am grateful to Dmitry Turaev for important conversations during the preparation of this paper, especially on the formulation of the arrayed standard blenders. I  thank Masayuki Asaoka, Xiaolong Li, Dmitrii Mints and Zhiyuan Zhang for useful discussions.  I also thank the anonymous referee for valuable comments and suggestions that improved the presentation of the paper. This work was supported by the Science Fund Program for Excellent
Young Scientists (Overseas) and the New Cornerstone Science Foundation.

\addcontentsline{toc}{section}{References}
\bibliographystyle{annotate}
\bibliography{references}

\end{document}